\documentclass[11pt]{amsart}
\usepackage{amscd,amsmath,amssymb,amsfonts}
\usepackage{mathabx}
\usepackage{color}
\usepackage{accents}
\usepackage[cmtip, all]{xy}
\usepackage{graphicx}
\usepackage{enumerate}
\theoremstyle{plain}
\newtheorem{thm}[equation]{Theorem}

\newtheorem{lem}[equation]{Lemma}

\newtheorem{prop}[equation]{Proposition}

\newtheorem{hyp}[equation]{Hypothesis}

\newtheorem{conj}[equation]{Conjecture}

\newtheorem{df}{Definition}[section]

\theoremstyle{definition}

\newtheorem{rem}{Remark}[section]

\newtheorem{rmks}[equation]{Remarks}
\newenvironment{dem}{\paragraph{Proof}}

\numberwithin{equation}{section}
\numberwithin{subsubsection}{subsection}

\newcommand{\isoarrow}{{~\overset\sim\longrightarrow}}

\bigskip

\newcommand{\N}{\mathbb{N}}
\newcommand{\Z}{\mathbb{Z}}
\newcommand{\Q}{\mathbb{Q}}
\newcommand{\C}{\mathbb{C}}
\newcommand{\R}{\mathbb{R}}

\newcommand{\AF}{\mathbb{A}_{F}}

\newcommand{\ZZ}{{\mathbb Z}}

\newcommand{\ra}{{~\rightarrow~}}

\newcommand{\QQ}{{\mathbb Q}}

\newcommand{\Qbar}{{\overline{\mathbb Q}}}

\newcommand{\RR}{{\mathbb R}}

\newcommand{\ad}{{{\mathbf A}}}

\newcommand{\CC}{{\mathbb C}}

\begin{document}

\author{Michael Harris}
\thanks{M.H.'s research received funding from the European Research Council under the European Community's Seventh Framework Programme (FP7/2007-2013) / ERC Grant agreement no. 290766 (AAMOT).  M.H. was partially supported by NSF Grant DMS-1404769.
}

\author{Jie LIN}


\thanks{J.L. was supported by the allocation de l'ENS and also by the European Research Council under the European Community's Seventh Framework Programme (FP7/2007-2013) / ERC Grant agreement no. 290766 (AAMOT)}

\address{Michael Harris\\
Department of Mathematics, Columbia University, New York, NY  10027, USA}
 \email{harris@math.columbia.edu}

\address{Jie Lin\\
Institut des Hautes \'{E}tudes Scientifiques, 91440 Bures-sur-Yvette}
 \email{jie.lin@imj-prg.fr}

\title [Period relations and special values of Rankin-Selberg $L$-functions]
{Period relations and special values of Rankin-Selberg $L$-functions}

\dedicatory{To Roger Howe, in friendship and admiration}

\maketitle

\begin{abstract}
This is a survey of recent work on values of Rankin-Selberg $L$-functions of pairs of cohomological automorphic representations that are {\it critical} in Deligne's sense.   The base field is assumed to be a CM field.  Deligne's conjecture is stated in the language of motives over $\QQ$, and express the critical values, up to rational factors, as determinants of certain periods of algebraic differentials on a projective algebraic variety over homology classes.   The results that can be proved by automorphic methods express certain critical values as (twisted) period integrals of automorphic forms.  Using Langlands functoriality between cohomological automorphic representations of unitary groups, which can be identified with the de Rham cohomology of Shimura varieties, and cohomological automorphic representations of $GL(n)$, the automorphic periods can be interpreted as motivic periods.  We report on recent results of the two authors, of the first-named author with Grobner, and of Guerberoff.
\end{abstract}

\section{Introduction}

Let $M$ be a motive of rank $n$ over $\QQ$,  which can be identified with a compatible
family of $\ell$-adic Galois representations $\rho_{\ell,M}$ of rank $n$.    Then  we define the $L$-function $L(s,M) = \prod_p L_p(s,M)$
   where for all $p$ for which $\rho_{\ell,M}$ is unramified
$$L_p(s,M) = [\det(1 - \rho_{\ell,M}(Frob_p)T)^{-1}]_{T = q^{-s}},$$ 
with $Frob_p$ any (geometric) Frobenius element in a decomposition group at $p$ inside $Gal(\Qbar/\QQ)$.   

Let $s_0 \in \ZZ$ be a {\it critical value} of $L(s,M)$, in Deligne's sense; we recall the definition below. We state 
(a crude version of) Deligne's conjecture regarding the value at $s_0$ of the $L$-function:    
$$L(s_0,M) \sim c^+(s_0,M)$$
where $c^+(s_0,M)$ is a certain explicitly defined determinant of periods of differential forms on (a smooth projective variety containing a realization of) $M(s_0)$, the twist of $M$ by the Tate motive $\QQ(s_0)$
and  $\sim $ means ``up to $\QQ$-multiples".   More precise versions of the conjecture are stated below for motives with coefficients in a number field $E$, and for motives over $\QQ$ obtained by restriction of scalars from a CM field $F$.

The theory of automorphic forms on reductive groups provides a large supply of $L$-functions, conjecturally including all the motivic $L$-functions $L(s,M)$.  Moreover, while the above definition of $L(s,M)$ suggests no obvious relation to differential forms on algebraic varieties, the special values of automorphic $L$-functions are often expressed as integrals of differential forms
over locally symmetric varieties.   Thus, practically all results on special values of motivic $L$-functions are in fact theorems about special values of automorphic $L$-functions, which can be identified with the $L$-functions of Galois representations obtained more or less directly from the theory of Shimura varieties.  

\subsection*{Examples}  In what follows, $F$ is a number field, and $GL(n)$, $GL(n-1)$ denote the algebraic groups $R_{F/\QQ}GL(n)_{F}$,  $R_{F/\QQ}GL(n-1)_{F}$ over $\QQ$.

\begin{itemize}
\item[(1)] It has been known for a long time that critical values of Rankin-Selberg $L$-functions of cohomological automorphic representations of $GL(n) \times GL(n-1)$ can often be written as cup products of differential forms on the associated adelic locally symmetric spaces.  Recent results on these lines are contained in \cite{Mah,rag,GH}.  The important paper of Binyong Sun \cite{sun}, which proves the non-vanishing of the relevant archimedean zeta integrals, shows that the cup product expressions can be used effectively to relate the critical values to natural period invariants obtained by comparing the rational structures defined by Whittaker models to those defined cohomologically.  There is no obvious relation, however, between these {\it Whittaker periods} and the {\it motivic periods} that enter into the computation of Deligne's period invariant $c^+$. 

\item[(2)]\label{cohh}  Suppose $F$ is a CM field, with maximal totally real field $F^+$, and let $\sigma \in Gal(F/F^+)$ be the non-trivial element.   A cuspidal cohomological automorphic representation $\Pi$ of $GL(n)$ that satisfies the hypothesis $\Pi^{\vee} \isoarrow \Pi^c$ gives rise by stable descent to an $L$-packet of cohomological automorphic representations $\{\pi\}_V$ of the group $U = R_{F^+/\QQ}U(V)$, where $V$ is an $n$-dimensional hermitian space over $F$, provided $\Pi$ is sufficiently ramified at the finite places at which $U$ is not quasi-split.  The members of $\{\pi\}_V$ define coherent cohomology classes on the Shimura variety $Sh(V)$ attached to $G$, and thus contribute to the Hodge-de Rham realizations of motives whose $L$-functions are related to the automorphic $L$-function $L(s,\Pi)$.

\item[(3)]  By applying the descent described in (2), the Rankin-Selberg $L$-function of $GL(n) \times GL(1)$ can be identified with the
{\it standard $L$-function} of the unitary group $U$,  whose integral representation by means of the doubling method identifies its critical values with periods of the coherent cohomology classes of $Sh(V)$.  This was carried out when $F^+ = \QQ$ in \cite{H97} and has recently been generalized by Guerberoff in \cite{gue} (see \ref{guer} below).

\item[(4)]  In many cases, by combining the methods of (1) with the results described in (3), one can express the Whittaker periods in terms of the periods of coherent cohomology classes, and thus with specific periods of motives realized in the cohomology of Shimura varieties -- we can call these {\it automorphic periods}.  As a result, the results of (1), specifically those of \cite{GH}, can be given motivic interpretations.  This has been generalized by Lin \cite{linthesis}, who has used these methods to obtain surprising {\it factorizations} of automorphic periods that are consistent with what is predicted by the Tate conjecture.

\end{itemize}

Some of these results are reviewed in the body of the paper.  We recall that Deligne's conjecture is stated for motives over $\QQ$ with coefficients in a field $E$; however, as Panchishkin observed \cite{panchishkin94}, motives over a totally real field $F$ with coefficients in $E$ can also be interpreted as motives over $\QQ$ with coefficients in $E\otimes F$, and Deligne's conjecture admits a refinement taking this into account.  The same is true when $F$ is a CM field.  Since this does not seem to have been treated in the literature, we include an appendix explaining the properties of $E\otimes F\otimes \CC$-modules that are relevant to Deligne's conjecture in our setting.

\subsection*{Acknowledgements}  We thank Harald Grobner and Lucio Guerberoff for helpful discussions regarding the topics of this survey, as well as for their many contributions to the contents of the paper.  We thank Li Ma for useful suggestions on the appendix of the paper. We also thank the editors for encouraging us to write the survey, and the referee for a very careful reading.  Finally, it is a pleasure and privilege to dedicate this paper to Roger Howe.

\section{Review of Deligne's Conjecture on critical values}

\subsection{Motives over general fields}
Let $M$ be a motive over a number field $F$ with coefficients in $E$. We denote by $\Sigma_{E}$ (resp. $\Sigma_{F}$) the set of embeddings of $E$ (resp. $F$) in $\C$. Tensor products without subscript are by default over $\Q$. The motive $M$ has several realizations as follows:

\begin{itemize}
\item Its de Rham realization $M_{DR}$ is a finitely generated free $E\otimes F$-module endowed with a Hodge filtration $M_{DR}\supset \cdots \supset F^{i}(M)\supset F^{i+1}(M)\supset \cdots$ where each $F^{i}(M)$ is a sub-$E\otimes F$-module of $M_{DR}$. We remark that $F^{i}(M)$ is in general not a free $E\otimes F$-module.

\item For each $\sigma\in \Sigma_{F}$, its Betti realization over $\sigma$ is a finite dimensional $E$-vector space $M_{\sigma}$. Moreover, we have a Hodge decomposition as $E\otimes \C$-modules:
\begin{equation}
M_{\sigma}\otimes \C=\bigoplus\limits_{p,q\in \Z}M^{p,q}_{\sigma}.\end{equation}

\item For each $\lambda$, a finite place of $E$, the $\lambda$-adic realization of $M$ is an $E_{\lambda}$-vector space endowed with an action of $Gal(\overline{F}/F)$ where $\overline{F}$ is an algebraic closure of $F$. The family $(M_{\lambda})_{\lambda}$ forms a compatible system of $\lambda$-adic representations. In particular, the $L$-function of $M$ can be defined as usual.
\end{itemize}

We refer to \cite{panchishkin94} for more discussions on motivic $L$-functions and $p$-adic motivic $L$-functions. Our goal here is to describe the Deligne conjecture and reformulate the Deligne period. For this purpose, we mainly focus on the de Rham realization and the Betti realizations.

The comparison isomorphism relates these realizations. More precisely, for each $\sigma\in\Sigma_{F}$, we have an isomorphism of free $E\otimes \C$-modules:
\begin{equation}
I_{\infty}: M_{\sigma}\otimes \C \xrightarrow{\sim} M_{DR}\otimes_{\sigma} \C.
 \end{equation}
 
This isomorphism is compatible with the Hodge structures in the sense that:
\begin{equation}\label{equation Hodge 1}
I_{\infty,\sigma}(\bigoplus\limits_{p\geq i} M^{p,q}_{\sigma})=F^{i}(M)\otimes_{\sigma} \C.
\end{equation}

From the isomorphisms above, we see that $dim_{E}M_{\sigma}=rank_{E\otimes F}M_{DR}$. We call this number the \textbf{rank} of $M$ and denote it by $\text{rank}(M)$.

For $w$ an integer, we say $M$ is \textbf{pure of weight} $w$ if $M^{p,q}_{\sigma}=0$ for any $\sigma$ and any $p,q$ such that $p+q\neq w$. 

Since $M_{\sigma}^{p,q}$ is an $E\otimes \C$-module, we can decompose it as $\bigoplus_{\tau\in \Sigma_{E}}M_{\sigma}^{p,q}(\tau)$ where the action of $E$ on $M_{\sigma}^{p,q}(\tau)$ is given by scalar multiplication via $\tau$.

We say $M$ is \textbf{regular} if $dim_{\C}M_{\sigma}^{p,q}(\tau) \leq 1$ for all $\sigma\in \Sigma_{F}$, $\tau\in \Sigma_{E}$ and $p,q\in \Z$.

We define the \textbf{Hodge type} of $M$ at $(\tau,\sigma)$ as the set $T(M_{\sigma})(\tau)$ consisting of pairs $(p,q)$ such that $M_{\sigma}^{p,q}(\tau)\neq 0$.

The \textbf{infinite Frobenius} at $\sigma$ is an $E$-linear isomorphism $F_{\infty,\sigma}: M_{\sigma} \rightarrow M_{\overline{\sigma}}$, which satisfies $F_{\infty,\overline{\sigma}}\circ F_{\infty,\sigma}=Id$. We may extend it to an $E\otimes \C$-linear isomorphism $M_{\sigma} \otimes \C \rightarrow M_{\overline{\sigma}}\otimes \C$. It sends $M_{\sigma}^{p,q}$ to $M_{\overline{\sigma}}^{q,p}$. Therefore if $(p,q)$ is contained in $T(M_{\sigma})(\tau)$, then $(q,p)$ is contained in $T(M_{\overline{\sigma}})(\tau)$.

We write $n$ for the rank of $M$. If $M$ is regular, then each $T(M_{\sigma})(\tau)$ contains exactly $n$ elements. If moreover $M$ is pure of weight $w$, then we can write the Hodge type of $M$ at $(\tau,\sigma)$ as $(p_{i}(\tau,\sigma),q_{i}(\tau,\sigma))_{1\leq i\leq n}$ with $p_{1}(\tau,\sigma) > p_{2}(\tau,\sigma)>\cdots >p_{n}(\tau,\sigma)$ and $q_{i}(\tau,\sigma)=w-p_{i}(\tau,\sigma)$.


\subsection{The Deligne conjecture}
The Deligne conjecture relates critical values of $L$-functions with the Deligne periods for motives over $\Q$. We will give a brief introduction here. We refer the reader to \cite{deligne79} for the details.

Let $\mathcal{M}$ be a motive over $\Q$ with coefficients in a number field $E$, pure of weight $w$. There is only one embedding of $\Q$ in $\C$. We write $\mathcal{M}_{B}$ for the Betti realization with respect to this unique embedding. It is an $E$-vector space and has a Hodge decomposition $\mathcal{M}_{B}\otimes \C= \bigoplus\limits_{p+q=\omega}\mathcal{M}^{p,q}$. Moreover, the infinite Frobenius $F_{\infty}$ exchanges $\mathcal{M}^{p,q}$ and $\mathcal{M}^{q,p}$.

The de Rham realization $\mathcal{M}_{DR}$ is endowed with a Hodge filtration of $E$-vector spaces.

The comparison isomorphism $I_{\infty}$ gives an isomorphism between $\mathcal{M}_{B}\otimes \C$ and $\mathcal{M}_{DR}\otimes \C$ as $E\otimes \C$-modules.

\begin{df}
We fix any $E$-bases of $\mathcal{M}_{B}$ (resp. $\mathcal{M}_{DR}$) and extend it to an $E\otimes \C$-basis of $\mathcal{M}_{B}\otimes \C$ (resp. $\mathcal{M}_{DR}\otimes \C$).

We define the \textbf{determinant period} with respect to the fixed bases $\delta(\mathcal{M})$ to be the determinant of the comparison isomorphism with respect to the fixed bases. It is an element in $(E\otimes \C)^{\times}$. Its image in $(E\otimes \C)^{\times}/E^{\times}$ does not depend on the choice of bases.

\end{df}

The comparison isomorphism induces an isomorphism between $\bigoplus\limits_{p\geq i}\mathcal{M}^{p,w-p}$ and $F^{i}\mathcal{M}\otimes \C$. This implies that $\mathcal{M}^{p,w-p}\cong (F^{i}\mathcal{M}/F^{i+1}\mathcal{M})\otimes \C$ is a free $E\otimes \C$-module. 

For simplicity, we assume that the middle stage of the filtration $\mathcal{M}^{w/2,w/2}=0$ if $w$ is even.

We write $\mathcal{M}^{+}$ (resp. $\mathcal{M}^{-}$) for $(\mathcal{M}_{B})^{F_{\infty}}$ (resp. $(\mathcal{M}_{B})^{-F_{\infty}}$), the subset of the fixed points of $F_{\infty}$ (resp. $-F_{\infty}$). Since $F_{\infty}$ is $E$-linear, it is an $E$-vector space. 

For $x\in \R$, we write $[x]$ for the largest integer which is no bigger than $x$. We write $F^{+}\mathcal{M}=F^{-}\mathcal{M}$ for $F^{[(w+1)/2]}\mathcal{M}$. The comparison isomorphism induces an $E\otimes \C$-linear map:
\begin{equation}
I_{\infty}^{\pm}:\mathcal{M}^{\pm}\otimes \C \hookrightarrow \mathcal{M}_{B}\otimes \C \cong \mathcal{M}_{DR}\otimes \C \rightarrow (\mathcal{M}_{DR}/F^{\pm}\mathcal{M})\otimes \C.
\end{equation}

One can show easily that $I_{\infty}^{\pm}\mathcal{M}$ is an isomorphism. 

\begin{df}
We fix any $E$-bases of $\mathcal{M}^{+}$ (resp. $\mathcal{M}_{DR}/F^{+}\mathcal{M}$) and extends it to an $E\otimes \C$-basis of $\mathcal{M}^{+}\otimes \C$ (resp. $(\mathcal{M}_{DR}/F^{+}\mathcal{M})\otimes \C$).

We define \textbf{the Deligne period} $c^{+}(\mathcal{M})$ with respect to the fixed bases to be the determinant of $I_{\infty}^{+}(\mathcal{M})$ with respect to the fixed bases. It is an element in $(E\otimes \C)^{\times}$. As before, its image in  $(E\otimes \C)^{\times}/E^{\times}$ does not depend on the choice of bases.

We may define \textbf{the Deligne period} $c^{-}(\mathcal{M})$ similarly.
\end{df}

\begin{df}
Let $A$ be an algebra (for example, $A$ is a field or the tensor product of two fields). Let $L$ be a subfield of $\C$.

Let $x$, $y$ be two elements in $A\otimes \C$. We write $x\sim_{A;L} y$ if there exists $e\in A\otimes L \subset A\otimes \C$ such that $x=ey$.

We write $\sim_{A}$ for $\sim_{A;\Q}$.
\end{df}

\begin{conj}(\textbf{the Deligne conjecture})
Let $m$ be a critical point for $\mathcal{M}$. We write $\epsilon$ for the sign of $(-1)^{m}$. We then have: 
\begin{equation}
L(m,\mathcal{M}) \sim_{E} (2\pi i)^{mn^{\epsilon}}c^{\epsilon}(\mathcal{M})
\end{equation}
where $n^{\epsilon}:=dim_{E}\mathcal{M}^{\epsilon}$.
\end{conj}

We refer to \cite{deligne79} for the definition of critical points. We also remark that in the case where $\mathcal{M}$ has no $(w/2,w/2)$-classes, we have $n^{+}=n^{-}=dim_{E}\mathcal{M}_{DR}/2$.

\subsection{Factorization of Deligne periods}

Let $F$ be a CM field of degree $d(F)$ over $\Q$. We write $F^+$ for the maximal totally real subfield of $F$.

Let $M$ be a rank $n$ motive over $F$ with coefficients in $E$, pure of weight $w$. We assume that $M^{w/2,w/2}=0$. Then for each $\sigma\in\Sigma_{F}$, we may define a $\sigma$-determinant period $\delta(M,\sigma)$ and a local Deligne period $c^{+}(M,\sigma)$ as follows.

Recall the comparison isomorphism:
\begin{equation}
I_{\infty,\sigma}: M_{\sigma}\otimes \C \xrightarrow{\sim} M_{DR}\otimes_{\sigma} \C.
\end{equation}

We fix an $E$-basis of $M_{\sigma}$ and extend it to an $E\otimes \C$-basis of $M_{\sigma}\otimes \C$. We fix an $E\otimes F$-basis of $M_{DR}$ and consider it as an $E\otimes \C$-basis of $M_{DR}\otimes_{\sigma} \C$. We define the \textbf{$\sigma$-determinant period} $\delta(M,\sigma)$ to be the determinant of $I_{\infty,\sigma}$ with respect to the fixed bases. It is an element in $(E\otimes \C)^{\times}$. As before, it depends on the choice of the bases, but its image in $(E\otimes\C)^{\times}/(E\otimes \sigma(F))^{\times}$ is independent of the choice.

\begin{rem}
We can use the same $E\otimes F$-basis for each $\sigma\in\Sigma_{F}$. Such a basis is called \textbf{covariant} as in Proposition $2.2$ of \cite{yoshidafactorization}. One can show that the image of the product $\prod\limits_{\sigma\in\Sigma_{F}}\delta(M,\sigma)$ in $(E\otimes \C)^{\times}/E^{\times}$ is independent of the choice of bases. Moreover, Proposition \ref{factorization theorem for motives} below holds only for a covariant basis.
\end{rem}

\bigskip

We now define the local Deligne periods $c^{+}(M,\sigma)$, extending the discussion in \cite{panchishkin94} to CM fields.  
 
The infinite Frobenius exchanges $M_{\sigma}$ and $M_{\overline{\sigma}}$. We fix $\Sigma$ a CM type of $F$, i.e. $\Sigma_{F}=\Sigma \bigsqcup \Sigma^{c}$. In particular, we know that $\Sigma$ has $d(F)/2$ elements.

For each $\sigma\in \Sigma$, we define $M_{\sigma}^{+}:=(M_{\sigma}\oplus M_{\overline{\sigma}})^{F_{\infty,\sigma}}$. It is an $E$-vector space of dimension $rank(M)$.

We write $F^{+}M$ for the $E\otimes F$-module $F^{[(w+1)/2]}M_{DR}$, and $M_{DR}^{+}$ for the $E\otimes F$-module $M_{DR}/F^{w/2}M_{DR}$. The comparison isomorphism induces an $E\otimes \C$-linear isomorphism:
\begin{equation}
I_{\infty,\sigma}^{+}: M_{\sigma}^{+}\otimes \C \xrightarrow{\sim} M_{DR}^{+}\otimes_{\sigma} \C \oplus M_{DR}^{+}\otimes_{\overline{\sigma}} \C= (M_{DR}^{+}\oplus M_{DR}^{+,c})\otimes_{\sigma}\C
\end{equation}
 where $M_{DR}^{+,c}$ is the same set as $M_{DR}^{+}$ endowed with the same action of $E$ and the complex conjugation of the action of $F$.

We know that $M_{DR}^{+}\oplus M_{DR}^{+,c}$ is a free $E\otimes F$-module by Lemma $2.1(3)$ of \cite{yoshidafactorization}. It can also be deduced easily from the fourth point of Proposition \ref{prop for decomposition} and the fact that $dim_{\C}(M_{DR}^{+}\oplus M_{DR}^{+,c})\otimes_{\tau \otimes\sigma} \C=dim M_{\sigma}^{+}\otimes_{\tau} \C=rank(M)=n$ does not depend on the choice of $\tau\in \Sigma_{E}$.

We define the \textbf{$\sigma$-Deligne period} $c^{+}(M,\sigma)$ to be determinant of $I_{\infty,\sigma}^{+}$ with respect to any fixed $E$-basis of $M_{\sigma}^{+}$ and a fixed $E\otimes F$-basis of $M_{DR}^{+}\oplus M_{DR}^{+,c}$. It is an element in $(E\otimes \C)^{\times}$. Its image in $(E\otimes \C)^{\times}/(E\otimes \sigma(F))^{\times}$ does not depend on the choice of the bases.

\begin{rem}
We may define $M_{\sigma}^{-}$, $M_{DR}^{-}$ and $c^{-}(M,\sigma)$ similarly. Recall that $M_{DR}^{-}=M_{DR}^{+}$ since $M$ has no $(w/2,w/2)$-classes. For such a motive, we have at complex places $\sigma$ that:
 \begin{equation}\label{relation between pm}
 c^{+}(M,\sigma)\sim_{E\otimes \sigma(F)}  c^{-}(M,\sigma).
 \end{equation}
  
 More precisely, let $\{e_{1},e_{2},\cdots,e_{n}\}$ be any $E$-basis of $M_{\sigma}$. Then $\{e_{i}+F_{\infty}e_{i}\}_{1\leq i\leq n}$ is an $E$-basis of $M_{\sigma}^{+}$ and $\{e_{i}-F_{\infty}e_{i}\}_{1\leq i\leq n}$ is an $E$-basis of $M_{\sigma}^{-}$. If we use these bases to calculate the $\sigma$-Deligne period, we will get:
 \begin{equation}\label{relation between pm}
 c^{+}(M,\sigma)= e_{\sigma} c^{-}(M,\sigma).
 \end{equation}
 Here $e_{\sigma}$ is an element in $E\otimes \C$ such that for each $\tau\in\Sigma_{E}$, $e_{\sigma}(\tau)=(-1)^{n_{\sigma}(\tau)}$ where 
 $n_{\sigma}(\tau)=dim_{\C}M_{DR}^{+,c}\otimes_{\sigma,\tau}\C$.

 For the proof, it suffices to consider the following commutative diagram:
$$\begin{CD}
M_{\sigma}^{+}\otimes\C @>I_{\infty,\sigma}^{+}>>  &M_{DR}^{+}\otimes_{\sigma} \C \oplus M_{DR}^{+,c}\otimes_{\sigma} \C \\
@VVV   &     @VVV \\
M_{\sigma}^{-}\otimes\C @>I_{\infty,\sigma}^{-}>>        &M_{DR}^{+}\otimes_{\sigma}\C\oplus M_{DR}^{+,c}\otimes_{\sigma} \C
\end{CD}$$
where the left vertical arrow is the $E\otimes \C$-linear isomorphism sending $e_{i}+F_{\infty}e_{i}$ to $e_{i}-F_{\infty}e_{i}$ and the right vertical arrow is the map $(Id,-Id)$.

It remains to show that $e_{\sigma}\in E\otimes \sigma(F)$. This is equivalent to show that for $g\in Gal(\bar{\Q}/\sigma(F))$ and $\tau\in\Sigma_{E}$, we have $e_{\sigma}(g\tau)=e_{\sigma}(\tau)$.

In fact, by Remark \ref{index is Galois invariant}, we know $n_{\sigma}(g\tau)=n_{g^{-1}\sigma}(\tau)=n_{\sigma}(\tau)$ where the last equality is due to the fact that $g^{-1}\sigma=\sigma$. Hence we have $e_{\sigma}(g\tau)=e_{\sigma}(\tau)$ as expected.
\end{rem}

\bigskip

We now consider the determinant period and the Deligne period for $\mathcal{M}=Res_{F/\Q}M$. We have: $\mathcal{M}_{DR}=M_{DR}$ and $\mathcal{M}_{B}=\bigoplus\limits_{\sigma\in\Sigma_{F}}M_{\sigma}$ as vector spaces over $E$.

\begin{prop}\label{factorization theorem for motives}
Let $\{w_{1},\cdots,w_{n}\}$ be an $E\otimes F$-basis of $M_{DR}$ such that the image of $\{(w_{i},w_{n+1-i})\}_{1\leq i\leq n}$ in $M_{DR}^{+}\oplus M_{DR}^{+,c}$ forms an $E\otimes F$-basis of $M_{DR}^{+}\oplus M_{DR}^{+,c}$. We use the family $\{w_{1},\cdots,w_{n}\}$ and the image of the family $\{(w_{i},w_{n+1-i})\}_{1\leq i\leq n}$ to define $\sigma$-determinant periods and $\sigma$-Deligne periods respectively for any $\sigma$.

Let $\alpha\in F$ be a purely imaginary element, i.e., $\overline{\alpha}=-\alpha$ where $\overline{\alpha}$ refers to the complex conjugation of $\alpha$ in the CM field $F$. The following factorizations of periods hold at the same time:
\begin{eqnarray}
\delta(\mathcal{M})\sim_{E} (D_{F}^{1/2})^{n}\prod\limits_{\sigma\in\Sigma_{F}} \delta(M,\sigma);\\
c^{+}(\mathcal{M})\sim_{E} (\prod\limits_{\sigma\in\Sigma}\sigma(\alpha))^{[n/2]}(D_{F^{+}}^{1/2})^{n}\prod\limits_{\sigma\in\Sigma} c^{+}(M,\sigma);\\
c^{-}(\mathcal{M})\sim_{E} (\prod\limits_{\sigma\in\Sigma}\sigma(\alpha))^{[n/2]}(D_{F^{+}}^{1/2})^{n}\prod\limits_{\sigma\in\Sigma} c^{-}(M,\sigma).
\end{eqnarray}
Here $D_{F}^{1/2}$ (resp. $D_{F^{+}}^{1/2}$) is the square root of the absolute discriminant of $F$ (resp. $F^{+}$). We identify it with $1\otimes D_{F}^{1/2}\in E\otimes \C$ (resp. $1\otimes D_{F^{+}}^{1/2}\in E\otimes \C$). Recall that $[n/2]$ is the largest integer no bigger than $n/2$.
\end{prop}

\begin{dem}
The first equation is proved in Proposition $2.2$ of \cite{yoshidafactorization}. The second one can be proved by similar argument. We now give the details.

The Deligne period $c^{+}(\mathcal{M})$ is the determinant of the composition of $\prod\limits_{\sigma\in\Sigma}I_{\infty,\sigma}^{+}$ and the following isomorphism of $E\otimes\C$-modules:
\begin{equation}\label{discriminant equation 1}
f: \bigoplus\limits_{\sigma\in\Sigma}(M_{DR}^{+}\oplus M_{DR}^{+,c})\otimes_{\sigma}\C \xrightarrow{\sim} M_{DR}^{+}\otimes \C.
\end{equation}

It remains to show that the determinant of the above isomorphism is equivalent to $ (\prod\limits_{\sigma\in\Sigma}\sigma(\alpha))^{[n/2]}D_{F^{+}}^{n/2}$ with respect to the fixed bases.

We know that $M_{DR}^{+}$ is free over $E\otimes F^{+}$. In fact, let $\widetilde{\sigma}$ be any element in $\Sigma_{F^{+}}$ and $\sigma$, $\overline{\sigma}$ be the places of $F$ over $\widetilde{\sigma}$. We know $dim_{\C} M_{DR}^{+}\otimes_{\tau\otimes \widetilde{\sigma}} \C= dim_{\C}(M_{DR}^{+}\otimes_{\tau \otimes\sigma} \C)+dim_{\C}(M_{DR}^{+}\otimes_{\tau \otimes\overline{\sigma}} \C) = n$ for any $\widetilde{\sigma}\in \Sigma_{F^{+}}$ and $\tau\in\Sigma_{E}$. Proposition \ref{prop for decomposition} then implies that $M_{DR}^{+}$ is a free $E\otimes F^{+}$-module of rank $n$.

We write $v_{i}$ for the image of $w_{i}$ in $M_{DR}^{+}$. 
We claim that the family
\begin{equation}\nonumber
\{v_{i}+v_{n+1-i}\}_{1\leq i\leq [(n+1)/2]} \cup\{ \alpha(v_{j}-v_{n+1-j})\}_{1\leq j\leq [n/2]}
\end{equation}
is an $E\otimes F^{+}$-basis of $M_{DR}^{+}$. By Lemma \ref{basislemma}, it is enough to prove that this family is linearly independent over $E\otimes F^{+}$.

We now prove this when $n=2m$ is even. In this case, if $\lambda_{i}$, $\mu_{i}$, $1\leq i\leq m$ are elements in $E\otimes F^{+}$ such that in $M_{DR}^{+}$ we have \begin{equation}\nonumber 
\sum\limits_{i=1}^{m}[\lambda_{i}(v_{i}+v_{2m+1-i})+\alpha\mu_{i}(v_{i}-v_{2m+1-i})]=0,
\end{equation}  hence
\begin{equation}\nonumber 
\sum\limits_{i=1}^{m}[(\lambda_{i}+\alpha\mu_{i})v_{i}+(\lambda_{i}-\alpha\mu_{i})v_{2m+1-i}]=0.
\end{equation} 
Recall that $\overline{\alpha}=-\alpha$. We know that in $M_{DR}^{+,c}$ \begin{equation}\nonumber \sum\limits_{i=1}^{m}[(\lambda_{i}-\alpha\mu_{i})v_{i}+(\lambda_{i}+\alpha\mu_{i})v_{2m+1-i}]=0.
\end{equation} We then deduce that \begin{equation}\nonumber 
\sum\limits_{1\leq i\leq m} (\lambda_{i}+\alpha\mu_{i})(v_{i},v_{2m+1-i})+\sum\limits_{1\leq i\leq m}(\lambda_{i}-\alpha\mu_{i})(v_{2m+1-i},v_{i})=0
\end{equation} in $M_{DR}^{+}\oplus M_{DR}^{+,c}$. Since $\{(v_{i},v_{2m+1-i})_{1\leq i\leq 2m}\}$ is an $E\otimes F$-basis of $M_{DR}^{+}\oplus M_{DR}^{+,c}$, we know $\lambda_{i}+\alpha\mu_{i}=\lambda_{i}-\alpha\mu_{i}$=0 for all $i$ and hence $\lambda_{i}=\mu_{i}=0$ as expected. The proof for odd $n$ is similar.

We can now take an $E$-basis of $M_{DR}^{+}$. Let $t_{1},\cdots, t_{d(F^{+})}$ be an integral basis of $F^{+}$ over $\Q$. Then $\{t_{k}(v_{i}+v_{n+1-i})\}_{1\leq i\leq [(n+1)/2],1\leq k\leq d(F^{+})} \cup\{ t_{k}\alpha(v_{j}-v_{n+1-j})\}_{1\leq j\leq [n/2],1\leq k\leq d(F^{+})}$ is an $E$-basis of $M_{DR}^{+}$.

The $E\otimes \C$-basis of $(M_{DR}^{+}\oplus M_{DR}^{+,c})\otimes_{\sigma}\C$ has been chosen as $\{(v_{i},v_{n+1-i})\otimes_{\sigma}1\}_{1\leq i\leq n}$. This basis is equivalent to \begin{equation}\nonumber 
\{(v_{i}+v_{n+1-i},v_{i}+v_{n+1-i})\otimes_{\sigma} 1\}_{1\leq i\leq [(n+1)/2]} \cup\{ (v_{j}-v_{n+1-j},-v_{j}+v_{n+1-j})\otimes_{\sigma} 1\}_{1\leq j\leq [n/2]}
\end{equation} by a rational transformation.

We observe that the $\sigma$-component of $f^{-1}(t_{k}(v_{i}+v_{n+1-i}))$ is $\sigma(t_{k})((v_{i}+v_{n+1-i},v_{i}+v_{n+1-i})\otimes_{\sigma} 1)$ for $1\leq i\leq [(n+1)/2]$ and $1\leq k\leq d(F^{+})$. The $\sigma$-component of $f^{-1}(t_{k}\alpha(v_{j}-v_{n+1-j}))$ is $\sigma(t_{k})(\sigma(\alpha)(v_{j}-v_{n+1-j}), \overline{\sigma}(\alpha)(v_{j}-v_{n+1-j}))\otimes_{\sigma}1 =\sigma({t_{k}})\sigma(\alpha)((v_{j}-v_{n+1-j},-v_{j}+v_{n+1-j})\otimes_{\sigma}1)$ for $1\leq j\leq [n/2]$ and $1\leq k\leq d(F^{+})$.

We then deduce that $det(f)^{-1}\sim_{E} (\prod\limits_{\sigma\in\Sigma}\sigma(\alpha))^{[n/2]}D_{F^{+}}^{n/2}.$
 
 Since $\prod\limits_{\sigma\in\Sigma}\sigma(\alpha)\times \prod\limits_{\sigma\in\Sigma}\overline{\sigma}(\alpha)\in \Q$ and hence $\prod\limits_{\sigma\in\Sigma}\sigma(\alpha)^{-1}\sim_{E}\prod\limits_{\sigma\in\Sigma}\overline{\sigma}(\alpha) \sim_{E} \prod\limits_{\sigma\in\Sigma}\sigma(\alpha)$ by the fact that $\overline{\alpha}=-\alpha$. We also have $D_{F^{+}}\in \Q$ and hence $D_{F^{+}}^{-1/2}\sim_{E} D_{F^{+}}^{1/2}$. We finally deduce that
 $det(f)\sim_{E} (\prod\limits_{\sigma\in\Sigma}\sigma(\alpha))^{[n/2]}D_{F^{+}}^{n/2}$ as expected.
\end{dem}

\begin{rem}
\begin{enumerate}
\item We shall see in the next section that such a basis always exists. Moreover, it can be good enough with respect to the Hodge decomposition.
\item If we change the condition to that the image of $\{(w_{i},w_{i})\}_{1\leq i\leq n}$ forms a basis in $M_{DR}^{+}\oplus M_{DR}^{+,c}$ and use this basis to calculate $\sigma$-Deligne periods, then we will have $c^{+}(\mathcal{M})\sim_{E} (D_{F^{+}}^{1/2})^{n}\prod\limits_{\sigma\in\Sigma} c^{+}(M,\sigma)$. Guerberoff has suggested a simple way to prove this. More precisely, we observe that $Res_{F^{+}/\Q}(Res_{F/F^{+}}(M))=Res_{F/F^{+}}M$. It is easy to see that $c^{+}(Res_{F/F^{+}}(M),\widetilde{\sigma})=c^{+}(M,\sigma)$. We may apply the factorization theorem for motives over totally real fields given in \cite{panchishkin94} and \cite{yoshidafactorization}. It remains to show that our basis is covariant for $Res_{F/F^{+}}(M)$. This is equivalent to saying that the image of the family $\{w_{i}\}_{1\leq i\leq n}$ is an $E\otimes F^{+}$-basis of $M_{DR}^{+}$, and this is not difficult to prove under the condition that the image of the family $\{(w_{i},w_{i})\}_{1\leq i\leq n}$ forms a basis in $M_{DR}^{+}\oplus M_{DR}^{+,c}$.
\item In general, for any basis, the relations $\delta(\mathcal{M})\sim_{E\otimes F^{gal}} \prod\limits_{\sigma\in\Sigma_{F}} \delta(M,\sigma)$ and $c^{+}(\mathcal{M})\sim_{E\otimes F^{gal}} \prod\limits_{\sigma\in\Sigma} c^{+}(M,\sigma)$ are always true where $F^{gal}$ is the Galois closure of $F$ over $\QQ$ (for example, see \cite{yoshidafactorization}). Changing the basis only affects the ratio of the two sides, which is an element of $F^{gal}$.
\end{enumerate}
\end{rem}

\section{Motivic and automorphic periods}

Notation is as in the previous section.  In this section we assume $M$ is regular.  Following \cite{linorsay}, we define period invariants generalizing those introduced in \cite{H97}.  Note that, in contrast to \cite{H97} and \cite{GH}, it is {\it not} assumed here that $M$ is polarized.

\subsection{Definitions} 

For $\sigma\in\Sigma_{F}$, we apply Proposition \ref{prop for decomposition} (2) to $F^{i}(M)$ and get:
\begin{equation}
F^{i}(M)(\alpha)\otimes_{\sigma} \C =\bigoplus\limits_{\tau \mid \alpha(\tau,\sigma)=\alpha} F^{i}(M)\otimes _{\tau\otimes \sigma} \C.
\end{equation}

On the other hand, equation (\ref{equation Hodge 1}) tells us that:
\begin{equation}
I_{\infty,\sigma}(\bigoplus\limits_{p\geq i} M^{p,q}_{\sigma}(\tau))=F^{i}(M)\otimes_{\tau\otimes\sigma} \C.
\end{equation}

We apply Proposition \ref{prop for decomposition} to $F^{i}(M)$ and get:
\begin{equation}\label{equation temp1}
I_{\infty,\sigma}(\bigoplus\limits_{p\geq i} M^{p,q}_{\sigma}(\tau))=F^{i}(M)(\alpha(\tau,\sigma))\otimes_{\tau\otimes\sigma} \C.
\end{equation}

Therefore,  \begin{equation}
dim_{\C}\bigoplus\limits_{p\geq i} M^{p,q}_{\sigma}(\tau)=dim_{L_{\alpha(\tau,\sigma)}}F^{i}M(\alpha(\tau,\sigma)).
\end{equation} 
Hence $dim_{\C}M^{p,q}_{\sigma}(\tau)=dim_{L_{\alpha(\tau,\sigma)}}F^{p}M(\alpha(\tau,\sigma))-dim_{L_{\alpha(\tau,\sigma)}}F^{p+1}M(\alpha(\tau,\sigma))$. In particular, the Hodge type at $(\tau,\sigma)$ only depends on $\alpha(\tau,\sigma)$.

For each $\alpha\in \mathcal{A}$ and $1\leq i\leq n$, we may define $p_{i}(\alpha)=p_{i}(\tau,\sigma)$ for any $(\tau,\sigma)$ such that $\alpha(\tau,\sigma)=\alpha$. We define $q_{i}(\alpha)=w-p_{i}(\alpha)$.

\bigskip

We may rewrite equation (\ref{equation temp1}) as:
\begin{equation}\label{temp equation 2}
I_{\infty,\sigma}(\bigoplus\limits_{j\leq i} M^{p_{j}(\alpha),q_{j}(\alpha)}_{\sigma}(\tau))=F^{p_{i}}(M)(\alpha)\otimes_{\tau\otimes\sigma} \C.
\end{equation}

\bigskip

In the rest of this chapter, we fix for every $\alpha\in\mathcal{A}$ and for every $i\in [1,n]$ an element $\widetilde{\omega_{i}}(\alpha)$ in $ F^{p_{i}(\alpha)}M(\alpha)\backslash F^{p_{i-1}(\alpha)}M(\alpha)$. Here we set $p_{0}=+\infty$ and $F^{p_{0}}(M)=\{0\}$.

Replacing $i$ by $i-1$ in equation (\ref{temp equation 2}), we get:\begin{equation}
I_{\infty,\sigma}(\bigoplus\limits_{j\leq i-1} M^{p_{j}(\alpha),q_{j}(\alpha)}_{\sigma}(\tau))=F^{p_{i-1}}(M)(\alpha)\otimes_{\tau\otimes\sigma} \C.
\end{equation}

We deduce that there exists $\epsilon_{i-1,\sigma}(\tau) \in F^{p_{i-1}}(M)(\alpha)\otimes_{\tau\otimes\sigma} \C$ such that $\omega_{i,\sigma}(\tau):=I_{\infty,\sigma}^{-1}(\widetilde{\omega_{i}}(\alpha)\otimes_{\tau\otimes \sigma}1-\epsilon_{i-1,\sigma}(\tau))$ is a non zero element in $M^{p_{i}(\alpha),q_{i}(\alpha)}_{\sigma}(\tau)$ where $\alpha=\alpha(\tau,\sigma)$. In particular, $\omega_{i,\sigma}(\tau)$ generates the one dimensional $\C$-vector space $M^{p_{i}(\alpha),q_{i}(\alpha)}_{\sigma}(\tau)=M^{p_{i}(\tau,\sigma),q_{i}(\tau,\sigma)}_{\sigma}(\tau)$.
 
For an integer $i\in[1,n]$, we write $i^{*}$ for $n+1-i$. The infinite Frobenius $F_{\infty,\sigma}$ maps $\omega_{i,\sigma}(\tau)$ to an element in $M^{p_{i^{*}}(\overline{\sigma},\tau),q_{i^{*}}(\overline{\sigma},\tau)}_{\sigma}(\tau)$. Therefore, there exists $Q_{i,\sigma}(\tau)\in \C^{\times}$ such that \begin{equation}
F_{\infty,\sigma}(\omega_{i,\sigma}(\tau)) = Q_{i,\sigma}(\tau) \omega_{i^{*},\overline{\sigma}}(\tau).
\end{equation} 

\begin{df}
We define the \textbf{motivic period} $Q_{i}(M,\sigma)$ as the element $(Q_{i,\sigma}(\tau))_{\tau\in\Sigma_{E}}$ in $(E\otimes \C)^{\times}$. 

For $1\leq j\leq n$, we also define \begin{equation}\label{def Qleq}
Q^{(j)}(M,\sigma):=Q_{1}(M,\sigma)Q_{2}(M,\sigma)\cdots Q_{j}(M,\sigma)\delta(M,\sigma)(2\pi i)^{n(n-1)/2}
\end{equation} and $Q^{(0)}(M,\sigma)=\delta(M,\sigma)(2\pi i)^{n(n-1)/2}$.
\end{df}

\begin{rem}
This definition depends on the choice of the fixed bases $\{\widetilde{\omega_{i}}(\alpha)\}_{1\leq i\leq n}$, $\alpha\in\mathcal{A}$. A different choice will change $Q_{i,\sigma}$ by a factor in $(E\otimes \sigma(F))^{\times}$.
\end{rem}

\begin{lem} \label{conjugacy on motivic period}We have the following equations on the motivic periods:
\begin{enumerate}
\item 
$\delta(M^{c},\sigma)=[\prod\limits_{i=1}^{n}Q_{i}(M,\sigma)]\delta(M,\sigma).$
\item For any $\sigma\in\Sigma_{F}$ and $0\leq j\leq n$, we have
\begin{equation}Q^{(n-j)}(M^{c},\sigma) \sim_{E\otimes \sigma(F)} Q^{(j)}(M,\sigma).
\end{equation}
\end{enumerate}
\end{lem}

\begin{dem}
We first show that (1) implies (2).
In fact, we have by definition that
\begin{eqnarray}\nonumber
Q^{(n-j)}(M^{c},\sigma) \sim_{E\otimes \sigma(F)} Q_{1}(M^{c},\sigma)Q_{2}(M^{c},\sigma)\cdots Q_{n-j}(M^{c},\sigma)\delta(M^{c},\sigma)(2\pi i)^{n(n-1)/2}\\
\text{ and }Q^{(j)}(M,\sigma)\sim_{E\otimes \sigma(F)}Q_{1}(M,\sigma)Q_{2}(M,\sigma)\cdots Q_{j}(M,\sigma)\delta(M,\sigma)(2\pi i)^{n(n-1)/2}\nonumber
\end{eqnarray}
It is easy to see that $Q_{i}(M^{c},\sigma) \sim_{E\otimes \sigma(F)} Q_{n+1-i}(M,\sigma)^{-1}$. Hence it remains to show (1).

In fact,  we have a commutative diagram:
\begin{equation}
\xymatrix{
M_{\sigma}\otimes\C \ar[d]^{F_{\infty}} \ar[r]^{I_{\infty,\sigma}} &M_{DR}\otimes_{\sigma}\C \ar[d]^{I_{\infty,\overline{\sigma}}F_{\infty,\sigma}I_{\infty,\sigma}^{-1}}\\
M_{\overline{\sigma}}\otimes\C \ar[r]^{I_{\infty,\overline{\sigma}}}          &M_{DR}\otimes_{\overline{\sigma}}\C}.
\end{equation}
Let $\omega_{i,\sigma}$ denote $\sum\limits_{\tau\in\Sigma_{E}}\omega_{i,\sigma}(\tau)$. By Proposition \ref{basis component}, the family $\{\omega_{i,\sigma}\}_{1\leq i\leq n}$ forms an $E\otimes \C$ basis of $M_{\sigma}\otimes \C$. Therefore, the family $\{I_{\infty}(\omega_{i,\sigma})\}_{1\leq i\leq n}$ forms a basis of  $M_{DR}\otimes_{\sigma}\C$. This basis is not rational, but can be transformed to a rational basis by a unipotent matrix. Hence it can be used to calculate the determinant period.

Similarly, we may use $\{I_{\infty}(\omega_{i,\overline{\sigma}})\}_{1\leq i\leq n}$ as the basis of $M_{DR}\otimes_{\overline{\sigma}}\C$.  Since $F_{\infty,\sigma}(\omega_{i,\sigma}) = Q_{i}(M,\sigma) \omega_{i^{*},\overline{\sigma}}$, the determinant of the right vertical arrow equals $\prod\limits_{i=1}^{n}Q_{i}(M,\sigma)$. The lemma then follows from the $E$-rationality of $F_{\infty}$.
\end{dem}

\subsection{Deligne period for tensor product of motives}
Let $M$ and $M'$ be two regular motives over a CM field $F$ with coefficients in a number field $E$ pure of weight $w$ and $w'$ respectively. We write $n$ for the rank of $M$ and $n'$ for the rank of $M'$.

We have defined motivic periods for $M$ and $M'$ in the previous sections. For each $\sigma\in\Sigma_{F}$, we shall calculate the local Deligne period for $M\otimes M'$ in terms of motivic periods in this section. We keep the notation of the last section.

We first construct an $E\otimes F$-basis of $M_{DR}$ which is good enough with respect to the Hodge filtration.

For each $\alpha\in\mathcal{A}$ and each $i\in [1,n]$, we define $\widehat{\omega_{i}}(\alpha):=\widetilde{\omega_{i}}(\alpha) -\sum\limits_{\alpha(\tau,\sigma)=\tau}\epsilon_{i-1,\sigma}(\tau)$. 
We consider $\{\widehat{\omega_{i}}(\alpha)\}_{1\leq i\leq n}$ as a family of vectors in $M(\alpha)\otimes\C$. Recall that $\{\widetilde{\omega_{i}}(\alpha)\}_{1\leq i\leq n}$ is an $L_{\alpha}\otimes \C$-basis of $M(\alpha)\otimes \C$. We claim that $\{\widehat{\omega_{i}}(\alpha)\}_{1\leq i\leq n}$ is also an $L_{\alpha}\otimes \C$-basis of $M(\alpha)\otimes\C$. 

In fact, since $\epsilon_{i-1,\sigma}(\tau) \in F^{p_{i-1}}(M)(\alpha)\otimes_{\tau\otimes\sigma} \C$, we know that the sum $\sum\limits_{\alpha(\tau,\sigma)=\tau}\epsilon_{i-1,\sigma}(\tau)$ is in $F^{p_{i-1}}(M)(\alpha)\otimes \C$. The $L_{\alpha}\otimes \C$-module $F^{p_{i-1}}(M)(\alpha)\otimes \C$ is generated by $\{\widetilde{\omega_{j}}(\alpha)\}_{1\leq j\leq i-1}$. Therefore, the family $\{\widehat{\omega_{i}}(\alpha)\}_{1\leq i\leq n}$ can be transformed to the $L_{\alpha}\otimes\C$-basis $\{\widetilde{\omega_{i}}(\alpha)\}_{1\leq i\leq n}$ by a unipotent matrix.\\

For each $1\leq i\leq n$, we define $\widetilde{\omega_{i}}$ (resp. $\widehat{\omega_{i}}$) to be the sum $\sum\limits_{\alpha\in\mathcal{A}}\widetilde{\omega_{i}}(\alpha)$ (resp. $\sum\limits_{\alpha\in\mathcal{A}}\widehat{\omega_{i}}(\alpha)$). The family $\{\widetilde{\omega_{i}}\}_{1\leq i\leq n}$ is a rational basis of $M_{DR}$ and hence can be considered as an $E\otimes F\otimes \C$-basis of $M_{DR}\otimes \C$.  The family $\{\widehat{\omega_{i}}\}_{1\leq i\leq n}$ is also an $E\otimes F \otimes \C$-basis of $M_{DR} \otimes \C$.  It can be transformed to the previous rational basis by a unipotent transformation. 

We remark that by the regularity property of $M$, the vector $\widehat{\omega_{i}}$ is unique up to multiplication by elements in $(E\otimes F)^{\times}$. \\

In order to apply Proposition \ref{factorization theorem for motives}, we first show that the basis $\{\widetilde{\omega_{i}}\}_{1\leq i\leq n}$ of $M_{DR}$ satisfies the condition there.

\begin{lem}
The image of the family $(\widetilde{\omega_{i}},\widetilde{\omega}_{n+1-i})_{1\leq i\leq n}$ in $M^{+}_{DR}\oplus M^{+,c}_{DR}$ forms an $E\otimes F$-basis.
\end{lem}

\begin{dem}
We write $v_{i}$ for the image of $\widetilde{\omega_{i}}$ in $M^{+}_{DR}=M_{DR}/F^{w/2}(M)$.

By Lemma \ref{basislemma}, it is enough to show that this family is $E\otimes F$-linearly independent. This is equivalent to saying that the family $\{(v_{i}\otimes_{\tau\otimes\sigma} 1, v_{i}\otimes_{\tau\otimes\overline{\sigma}}1)\}_{1\leq i\leq n}$ is $\C$-linearly independent in $M^{+}_{DR}\otimes_{\tau\otimes\sigma}\C \oplus M^{+}_{DR}\otimes_{\tau\otimes\overline{\sigma}}\C$.

Recall that by construction we have $v_{i}\otimes_{\tau\otimes\sigma} 1\in F^{p_{i}(\tau,\sigma)}M\otimes_{\tau\otimes\sigma}\C$. Moreover, $v_{n+1-i}\otimes_{\tau\otimes\overline{\sigma}} 1\in F^{p_{n+1-i}(\tau,\overline{\sigma})}M\otimes_{\tau\otimes\overline{\sigma}}\C=F^{w-p_{i}(\tau,\sigma)}M\otimes_{\tau\otimes\sigma}\C$.

Since $p_{i}(\tau,\sigma)\neq w/2$, we have either $p_{i}(\tau,\sigma)> w/2$, either $p_{i}(\tau,\sigma)< w/2$. In the first case, we have $F^{p_{i}(\tau,\sigma)}M\subset F^{w/2}(M)$ and hence $v_{i}\otimes_{\tau\otimes\sigma}1=0$. In the second case we have $v_{n+1-i}\otimes_{\tau\otimes\overline{\sigma}}1=0$.

We deduce that the family $\{(v_{i}\otimes_{\tau\otimes\sigma} 1, v_{i}\otimes_{\tau\otimes\overline{\sigma}}1)\}_{1\leq i\leq n}$ equals $\{(v_{i}\otimes_{\tau\otimes\sigma}1,0)\}_{p_{i}(\tau,\sigma)<w/2}\cup \{(0, v_{n+1-i}\otimes_{\tau\otimes\overline{\sigma}}1)_{p_{i}(\tau,\sigma)\}>w/2}$. The linear independence of this family is clear by the construction of $\widetilde{\omega}_{i}$.

\end{dem}

\bigskip

Therefore, we may apply Proposition \ref{factorization theorem for motives} to the basis $\{\widetilde{\omega_{i}}\}_{1\leq i\leq n}$ and calculate the Deligne period using this basis. Since the determinant of a unipotent matrix is always one, we can instead use the basis $\{\widehat{\omega_{i}}\}_{1\leq i\leq n}$. The idea was implicitly contained in \cite{H97} and discussed in detail in \cite{linorsay} and \cite{gue}.

We fix a basis for $M'_{DR}$ similarly. For any $\sigma\in\Sigma_{F}$, we fix any $E$-bases of $M_{\sigma}$ and $M'_{\sigma}$ and extend them to $E\otimes \C$-bases of $M_{\sigma}\otimes \C$ and $M'_{\sigma}\otimes \C$ respectively.

It remains to define the split index to state our main proposition.
\begin{df}
We write the Hodge type of $M$ at $(\tau,\sigma)$ as $(p_{i}(\tau,\sigma),w-p_{i}(\tau,\sigma))_{1\leq i\leq n}$ with $p_{1}(\tau,\sigma)>p_{2}(\tau,\sigma)>\cdots >p_{n}(\tau,\sigma)$. We write the Hodge type of $M'$ at $(\tau,\sigma)$ as $(r_{j}(\tau,\sigma),w'-r_{j}(\tau,\sigma))_{1\leq j\leq n'}$ with $r_{1}(\tau,\sigma)>r_{2}(\tau,\sigma)>\cdots >r_{n'}(\tau,\sigma)$.

We assume that $M\otimes M'$ has no $(\frac{w+w'}{2},\frac{w+w'}{2})$-classes. Then $p_{i}(\tau,\sigma)+r_{j}(\tau,\sigma)\neq \frac{w+w'}{2}$ for any $i$, $j$. Hence the sequence $-r_{n'}(\tau,\sigma)>-r_{n'-1}(\tau,\sigma)>\cdots>-r_{1}(\tau,\sigma)$ is split into $n+1$ parts by the numbers $p_{1}(\tau,\sigma)-\frac{w+w'}{2}>p_{2}(\tau,\sigma)-\frac{w+w'}{2}>\cdots>p_{n}(\tau,\sigma)-\frac{w+w'}{2}$. We denote the length of each part by $sp(i,M;M',\sigma)(\tau)$, $0\leq i\leq n$, and call them the \textbf{split indices} for the motivic pair.

For $0\leq i\leq n$, we write $sp(i,M;M',\sigma)$ for $(sp(i,M;M',\sigma)(\tau))_{\tau\in\Sigma_{E}}$ as an element in $\N^{\Sigma_{E}}$. We may define $sp(j,M';M,\sigma)$ for $1\leq j\leq n'$ similarly.
\end{df}

\begin{prop}\label{tensor product factorization}
We assume that the motive $M\otimes M'$ has no $(\frac{w+w'}{2},\frac{w+w'}{2})$-classes. 
We have the following equation for the Deligne period with respect to the above bases:
\begin{eqnarray}  &c^{+}(M\otimes M',\sigma) =
&\\ 
&(2\pi i)^{-\frac{nn'(n+n'-2)}{2}}\prod\limits_{j=0}^{n}(Q^{(j)}(M,\sigma))^{sp(j,M;M',\sigma)}\prod\limits_{k=0}^{n'}(Q^{(k)}(M',\sigma))^{sp(k,M';M,\sigma)}\nonumber&
\end{eqnarray}
Moreover, the factorization of $c^{+}(Res_{F/\Q}(M\otimes M'))$ in Proposition \ref{factorization theorem for motives} holds at the same time.
\end{prop}

Since the previous proposition concerns only one place, we may assume that the base field $F^+ = \QQ$. In this case, the previous proposition was proved in \cite{harrisadjoint, GH} when the motive is polarized, and was generalized to non-polarized motives in \cite{linthesis,linorsay}. 

\begin{rem}
If $M\otimes M'$ has non-trivial $(w/2,w/2)$-classes then there is no critical point (cf. $1.7$ of \cite{H97}).
\end{rem}

Recall by definition that the image of $\widehat{\omega_{i}}$ in $M_{DR}\otimes_{\sigma} \C$ is equal to $I_{\infty,\sigma}(\omega_{i,\sigma})$. The proof is the same as for motives over quadratic imaginary fields (cf. Proposition $1.1$ and Proposition $1.2$ of \cite{linorsay}). We only need to replace $\omega_{i}$ (resp. $\omega_{i}^{c}$) there by our $\omega_{i,\sigma}$ (resp. $\omega_{i,\overline{\sigma}}$). 

The key point of the proof is the relation $F_{\infty,\sigma}(\omega_{i,\sigma})= Q_{i,\sigma}\omega_{i^{*},\overline{\sigma}}$. The original ideas can be found in \cite{harrisadjoint} where the author provided a proof for self dual motives.

\subsection{The Deligne period for automorphic pairs}

Let $\Pi$ and $\Pi'$ be cuspidal cohomological automorphic representations of $GL_{n}(\AF)$ and $GL_{n'}(\AF)$ respectively.   We write $(z^{A_{\sigma,i}}\overline{z}^{A_{\overline{\sigma},i}})_{1\leq i\leq n}$ for the infinity type of $\Pi$ and $(z^{B_{\sigma,j}}\overline{z}^{B_{\overline{\sigma},j}})_{1\leq j\leq n'}$ for the infinity type of $\Pi'$ at $\sigma\in\Sigma$ respectively. The numbers $A_{\sigma,i}$ are in $ \Z+\frac{n-1}{2}$ and the numbers $B_{\sigma,j}$ are in $ \Z+\frac{n'-1}{2}$ for any $\sigma\in\Sigma_{F}$, $1\leq i\leq n$ and $1\leq j\leq n'$; they are written in strictly decreasing order:
$$A_{\sigma,i} > A_{\sigma,i+1}; ~~  B_{\sigma,j} > B_{\sigma,j+1}$$
for all $\sigma$ and $1 \leq i < n$, $1 \leq j < n'$.   In particular, $\Pi$ and $\Pi'$ are {\it regular}.    

We know that the sum $A_{\sigma,i}+A_{\overline{\sigma},i}=w(\Pi)$ does not depend on the choice of $i$ or $\sigma$. Moreover, the finite part of $\Pi$ is defined over a number field $E(\Pi)$. We define $w(\Pi')$ and $E(\Pi')$ similarly.


\begin{df}\label{thehypotheses}
\begin{enumerate}
\item We say $\Pi$ is \textbf{polarized} if $\Pi^{\vee}\cong \Pi^{c}$.
\item We say $\Pi$ is \textbf{sufficiently regular} if $\mid A_{\sigma,i}-A_{\sigma,i'}\mid$ is big enough\footnote{Different results require gaps of different sizes.  For example, in Theorem \ref{RSoverQ}, the gap has to be just slightly larger than the minimum; larger gaps are needed for some of the results of \cite{linthesis}}  for any $\sigma\in\Sigma_{F}$ and $i\neq i'$.
\item We say the pair $(\Pi,\Pi')$ is in \textbf{good position} if $n>n'$ and the numbers $-B_{j}$, $1\leq j\leq n'$ lie in different gaps between the numbers $-\frac{w(\Pi)+w(\Pi')}{2}+A_{i}$, $1\leq i\leq n$.
\end{enumerate}
\end{df}

We write $M(\Pi)$ and $M(\Pi')$ for the pure motives over $F$ (for absolute Hodge cycles) conjecturally attached to $\Pi$ and $\Pi'$ respectively. We know $M(\Pi)$ (resp. $M(\Pi')$) is pure of weight $w:=-w(\Pi)+n-1$ (resp. $w':=-w(\Pi')+n'-1$) with coefficients in $E(\Pi)$ (resp. $E(\Pi')$).  We will consider the Deligne conjecture when $Res_{F/\Q}(M(\Pi)\otimes M(\Pi'))$ has no $(\frac{w+w'}{2},\frac{w+w'}{2})$-classes.

\begin{df}
For $1\leq i\leq n$, $1\leq j\leq n'$ and $\sigma\in\Sigma$, we define $sp(i,\Pi;\Pi',\sigma):=sp(i,M(\Pi);M(\Pi'),\sigma)$ and $sp(j,\Pi';\Pi,\sigma):=sp(j,M(\Pi');M(\Pi),\sigma)$.
\end{df}

\begin{rem}
\begin{enumerate}
\item When we write down the infinity type of $\Pi$, we have implicitly fixed an embedding of the coefficient field $E$.  Changing the embedding induces a permutation on the infinity type. For example, let $\chi$ be an algebraic Hecke character with coefficients in a number field $E$. We fix $\tau: E\hookrightarrow \C$ an embedding of $E$. We can then write its infinity type as $(\sigma(z)^{a_{\sigma}})_{\sigma\in\Sigma_{F}}$. Let $\tau'$ be another embedding of $E$. We take $g\in Aut(\C)$ such that $\tau'=g\circ \tau$. Then the infinity type of $\chi$ with respect to the embedding $\tau'$ is $(\sigma(z)^{a_{g^{-1}\circ \sigma}})_{\sigma\in\Sigma_{F}}$ .
\item
The Hodge type of $M(\Pi)$ at $\sigma$ with respect to a fixed embedding of $E$ is $\{(-A_{\sigma,i}+\frac{n-1}{2},-A_{\overline{\sigma},i}+\frac{n-1}{2})\}$. Hence the integers $sp(i,\Pi;\Pi',\sigma)$ can be defined without assuming the existence of $M(\Pi)$.
\end{enumerate}
\end{rem}

By Proposition \ref{tensor product factorization}, we have 
\begin{eqnarray}  \label{local motivic period calculation}&c^{+}(M(\Pi)\otimes M'(\Pi'),\sigma) 
&\\ 
&=(2\pi i)^{\frac{-nn'(n+n'-2)}{2}}\prod\limits_{j=0}^{n}(Q^{(j)}(M(\Pi),\sigma))^{sp(j,\Pi;\Pi',\sigma)}\prod\limits_{k=0}^{n'}(Q^{(k)}(M(\Pi'),\sigma))^{sp(k,\Pi';\Pi,\sigma)}\nonumber&
\end{eqnarray}

If $n'=1$, it is easy to calculate the split index. 

\begin{prop}\label{n1motive}
Let $\Pi'=\chi$ be a regular algebraic Hecke character over $F$ of infinity type $(\sigma(z)^{a_{\sigma}})_{\sigma\in\Sigma_{F}}$.

For $\sigma$ in the CM type $\Sigma$, we write $I_{\sigma}=I_{\sigma}(\Pi,\chi)$ for the cardinal of the set $\{i\mid A_{\sigma,i}-A_{\overline{\sigma},i}+a_{\sigma}-a_{\overline{\sigma}}<0\}$.

If $M(\Pi)\otimes M(\chi)$ is critical, then \begin{eqnarray}
&c^{+}(M(\Pi)\otimes M(\chi),\sigma)\sim_{E(\Pi)\otimes E(\chi);\sigma(F)} \nonumber&\\
&(2\pi i)^{-\frac{n(n-1)}{2}} Q^{(I_{\sigma})}(M(\Pi),\sigma) Q^{(0)}(M(\chi),\sigma)^{n-I_{\sigma}}Q^{(1)}(M(\chi),\sigma)^{I_{\sigma}}\nonumber&
\end{eqnarray}
\end{prop}

For each $\sigma\in\Sigma_{F}$, the CM period $p(\chi,\sigma)$ is a complex number defined in the appendix of \cite{harrisappendix}. 

We define $\widecheck{\chi}=\chi^{c,-1}$ and $\widetilde{\chi}=\cfrac{\chi}{\chi^{c}}$. It is easy to see that $\widecheck{\widetilde{\chi}}=\widetilde{\chi}$.

The following lemma follows from Blasius's results on the Deligne conjecture \cite{blas}.

\begin{lem}(Comparison of CM periods and motivic periods for Hecke characters)\label{comparison n=1}
Let $\chi$ be a regular algebraic Hecke character over $F$. We have:
\begin{enumerate}
\item $\delta(M(\chi),\sigma) \sim_{E(\chi)} p(\widecheck{\chi^{c}},\sigma)$;
\item $Q_{1}(M(\chi),\sigma) \sim_{E(\chi)} \cfrac{p(\widecheck{\chi},\sigma)}{p(\widecheck{\chi^{c}},\sigma)} \sim_{E(\chi)}p( \cfrac{\widecheck{\chi}}{\widecheck{\chi^{c}}},\sigma)$.

Consequently, we have:
\begin{equation}\label{motivic period for Hecke characters}
Q^{(0)}(M(\chi),\sigma)\sim_{E(\chi)} p(\widecheck{\chi^{c}},\sigma) \text{ and }  Q^{(1)}(M(\chi),\sigma)\sim_{E(\chi)} p(\widecheck{\chi},\sigma).
\end{equation}
\end{enumerate}

\end{lem}
We refer to section $6.4$ of \cite{linthesis} for the proof of the lemma when $F^{+}=\QQ$. The same ideas should work for general CM fields.

\bigskip

Proposition \ref{n1motive} generalizes the similar expression proved in \cite{H97} when $F^+ = \QQ$, and can be compared to Guerberoff's expression in \cite{gue} for the Deligne period of $M_{0} \otimes R_{F/F^{+}}M(\widetilde{\psi})$;  here $M_{0}$ is a polarized motive over $F^{+}$ and $\psi$ is an algebraic Hecke character of $F$ of infinity type $(z^{-m_{\sigma}})_{\sigma\in\Sigma_{F}}$. For simplicity, we assume that $\psi$ is of weight $0$, i.e. $m_{\sigma}+m_{\overline{\sigma}}=0$ for all $\sigma_{0}\in\Sigma_{F}$. 

\begin{thm}\label{guermotive}[Guerberoff] If $M_{0}\otimes Res_{F/F^{+}}M(\widetilde{\psi})$ is critical, then for any $\sigma_{0}\in\Sigma_{F^{+}}$, we have
\[ c^{+}\left(M_{0} \otimes RM(\widetilde{\psi}),\sigma_{0}\right)\sim_{E\otimes \Q(\widetilde{\psi})\otimes\sigma(F)   }
\delta(M_{0},\sigma_{0})Q(\widetilde{\psi},\sigma_{0})\prod_{j\leq s_{\sigma}}Q_{j,\sigma_{0}}(M_{0}).\]
\end{thm}

We refer to \cite{gue} for an explanation of the parameters for the general case. 

If in addition to being polarized, $\Pi$ is isomorphic to $\Pi^{c}$, we expect that $\Pi$ is associated to a polarized motive $M_{0}$ over $F^{+}$.

In this case we have \begin{equation}\nonumber
c^{+}(M(\Pi)\otimes M(\widetilde{\psi}),\sigma)\sim_{E\otimes \Q(\widetilde{\psi})\otimes \sigma(F)} c^{+}(M_{0}\otimes Res_{F/F^{+}}M(\widetilde{\psi})).
\end{equation}

We claim that Guerberoff's result is compatible with Proposition \ref{n1motive}. 

In fact, we write $I_{\sigma}$ as the index in Proposition \ref{n1motive} for the pair $(\Pi,\widetilde{\psi})$. One can show that the index $s_{\sigma}$ in Theorem \ref{guermotive} is equal to $n-I_{\sigma}$. Hence the index $r_{\sigma}$ in \cite{gue} is equal to $I_{\sigma}$. As in section $5$ of \cite{gue}, we have 
\begin{equation}\nonumber
Q(\widetilde{\psi},\sigma_{0})\sim_{ \Q(\widetilde{\psi})\otimes \sigma(F)}  p(\widetilde{\psi},\sigma)^{r_{\sigma}-s_{\sigma}}=p(\widetilde{\psi},\sigma)^{2I_{\sigma}-n}.
\end{equation}

On the other hand, by equation (\ref{motivic period for Hecke characters}) we may simplify the last two terms in Proposition \ref{n1motive} for $\chi=\widetilde{\psi}$ as follows:
\begin{equation}\nonumber
Q^{(0)}(M(\widetilde{\psi}),\sigma)^{n-I_{\sigma}}Q^{(1)}(M(\widetilde{\psi}),\sigma)^{I_{\sigma}}\sim_{\Q(\widetilde{\psi})\otimes \sigma(F)} p(\widetilde{\psi},\sigma)^{2I_{\sigma}-n}.
\end{equation}

It remains to compare $(2\pi i)^{-\frac{n(n-1)}{2}} Q^{(I_{\sigma})}(M(\Pi),\sigma)$ with $\delta(M_{0},\sigma_{0})\prod_{j\leq s_{\sigma}}Q_{j,\sigma_{0}}(M_{0})$. It follows from the definitions that 
$$\delta(M_{0},\sigma_{0})\sim_{E\otimes\sigma(F)}\delta(M(\Pi),\sigma)$$
 and 
 $$Q_{j}(M_{0},\sigma_{0})\sim_{E\otimes \sigma(F)} Q_{j}(M_{0},\sigma)$$ 
 where $\sigma$ is the lifting of $\sigma_{0}$ in the CM type $\Sigma$ of $F$. Hence the term $\delta(M_{0},\sigma_{0})\prod_{j\leq s_{\sigma}}Q_{j,\sigma_{0}}(M_{0})$ is equivalent to $(2\pi i)^{-\frac{n(n-1)}{2}} Q^{(s_{\sigma})}(M(\Pi),\sigma)$. Since $\Pi\cong \Pi^{c}$, we have $M(\Pi) \cong M(\Pi)^{c}$. By the second part of Lemma \ref{conjugacy on motivic period} we have $Q^{(s_{\sigma})}(M(\Pi),\sigma)\sim_{E\otimes \sigma(F)} Q^{(n-s_{\sigma})}(M(\Pi)^{c},\sigma) \sim_{E\otimes \sigma(F)} Q^{(I_{\sigma})}(M(\Pi),\sigma)$ which completes the comparison.







\section{Results on critical values}

We can now restate the Deligne conjecture for $Res_{K/\Q}(M(\Pi)\otimes M(\Pi'))$.
\begin{conj}\label{deligneMMprime}
If $m\in \N+\frac{n+n'-2}{2} $ is critical for $\Pi\times \Pi'$ then the critical value at $s=m$ satisfies
\begin{eqnarray}
&L(m,\Pi\times \Pi') \sim_{E(\Pi)\otimes E(\Pi'); F^{gal}}&\\
 &(2\pi i)^{mnn'd(F^{+})}\prod\limits_{\sigma\in\Sigma} [\prod\limits_{j=0}^{n}(Q^{(j)}(M(\Pi),\sigma))^{sp(j,\Pi;\Pi',\sigma)}\prod\limits_{k=0}^{n'}(Q^{(k)}(M(\Pi'),\sigma))^{sp(k,\Pi';\Pi,\sigma)}].&\nonumber
\end{eqnarray}
\end{conj}

\begin{rem}Here we consider the relation up to $E(\Pi)\otimes E(\Pi')\otimes F^{gal}$. Hence we can ignore the terms $\prod\limits_{\sigma\in\Sigma}\sigma(\alpha)$ and $D_{F^{+}}^{1/2}$ in Proposition \ref{factorization theorem for motives}. Moreover, by equation (\ref{relation between pm}), we don't need to consider the sign in the original Deligne conjecture.
\end{rem} 
\begin{rem}
When we write down the relation $\sim_{E(\Pi)\otimes E(\Pi');F^{gal}}$, we have implicitly considered $L(m,\Pi\otimes\Pi')$ as an element of $E(\Pi)\otimes E(\Pi') \otimes \C$. In fact, we have a priori fixed embeddings of $E(\Pi)$ and $E(\Pi')$ in $\C$. For each $g_{1}\in Aut(\C)$ and $g_{2}\in Aut(\C)$, the value $L(m,\Pi^{g_{1}}\otimes\Pi'^{g_{2}})$ depends only on $g_{1}\mid_{E(\Pi)}$ and $g_{2}\mid_{E(\Pi')}$. Therefore, for any embeddings $\tau_{1}: E(\Pi) \hookrightarrow \C$ and $\tau_{2}: E(\Pi) \hookrightarrow \C$, we may define $L(m,\Pi^{\tau_{1}}\otimes\Pi'^{\tau_{2}})$ as $L(m,\Pi^{g_{1}}\otimes\Pi'^{g_{2}})$ where $g_{i}$ is any lift of $\tau_{i}$ in $Aut(\C)$ for $i=1,2$. It is clear that $(L(m,\Pi^{\tau_{1}}\otimes\Pi'^{\tau_{2}}))_{\tau_{1}: E(\Pi) \hookrightarrow \C, \tau_{2}: E(\Pi') \hookrightarrow \C}$ is an element of $E(\Pi)\otimes E(\Pi')\otimes\C$. We simply denote it by $L(m,\Pi\otimes\Pi')$.
\end{rem}

Some concrete results have been shown when the motivic periods $Q^{(j)}(M(\Pi),\sigma)$ are replaced by corresponding automorphic arithmetic periods. Here is a list of results obtained recently when $F$ is a CM field.  Notation ($P^{(a)}(\Pi)$, etc.) will be explained more precisely in subsequent sections.   In what follows, we let $K$ be a quadratic imaginary field, and assume $F=F^{+}K$ where $F^{+}$ is the maximal totally real subfield of $F$.

\subsection{Results when $F^+ = \QQ$}

We first assume $F^+ = \QQ$, where the results are easier to state.

\begin{thm}\label{RSoverQ}[Grobner-H., Lin]  We assume that both $\Pi$ and $\Pi'$ are cuspidal, cohomological, polarized and sufficiently regular. 

If $n \nequiv n' \mod 2$ and the pair $(\Pi,\Pi')$ is in good position (cf. Definition \ref{thehypotheses}) then for critical points $m>0$ we have that:
$$L(m,\Pi \times \Pi') \sim_{E(\Pi)\otimes E(\Pi'); F^{gal}} (2\pi i)^{mnn'}  \prod_{j = 0}^{n} P^{(j)}(\Pi)^{sp(j,\Pi;\Pi')}\prod_{k = 0}^{n'} P^{(k)}(\Pi')^{sp(k,\Pi';\Pi)},$$
where $P^{(k)}(\Pi) = P^{(I)}(\Pi)$, in the notation of Definition \ref{autoperiods}, with $I$ the singleton $k \in \{0,\dots, n\}$.
\end{thm}  






\subsection{Guerberoff's results on critical values of $L$-functions}

As noted above, Guerberoff has generalized the results of \cite{H97} on special values of $L$-functions of unitary groups to arbitrary CM fields.  The statement of Guerberoff's result is as in \cite{gue}; some of the notation has been changed.

\begin{thm}\label{guer}[Guerberoff] 
Let $V$ be a Hermitian space of dimension $n$ over $F$ respect to $F/F^{+}$. Assume that $V$ has signature $(r_{\sigma},s_{\sigma})$ at $\sigma \in \Sigma$ where $\Sigma$ is a CM type of $F$ as before.

Let $\pi$ be a cohomological holomorphic discrete series automorphic representation of the rational similitude group associated to $V$. Let $\psi$ be an algebraic Hecke character of $F$ of infinity type $(z^{-m_{\sigma}}))_{\sigma\in\Sigma_{F}}$.   We can parametrize the weight of $\pi$ by a tuple of integers $((a_{\sigma,1},\dots,a_{\sigma,n})_{\sigma\in\Sigma};a_{0})$ in a natural way.  

We assume that $\pi$ is polarized, i.e. $\pi^{\vee}\cong \pi\otimes ||\mu||^{2a_{0}}$.
 If $ m\in \Z+\frac{n-1}{2}$ is such that 
\begin{equation}\label{critical condition gue}
n<m+\frac{n-1}{2}  \leq\min\{a_{\sigma,r_{\sigma}}+s_{\sigma}+m_{\sigma}-m_{\overline\sigma},a_{\sigma,s_{\sigma}}+r_{\sigma}+m_{\overline\sigma}-m_{\sigma}\}_{\sigma\in\Sigma},
\end{equation}
then
\[ L^{S}\left(m,\pi\otimes\psi, St\right)\sim_{E(\pi)\otimes E(\psi);F^{gal}}(2\pi i)^{d(F^{+})mn-2a_{0}}P(\psi)Q_V(\pi).\]
\end{thm}

The term $P(\psi)$ is a certain expression involving periods of CM abelian varieties attached to the Hecke character $\psi$.  The term $Q_V(\pi)$ is the normalized Petersson norm of an arithmetic holomorphic vector in $\pi$.  

\begin{rmks}  (1)  The current version of Guerberoff's theorem only applies to polarized representations, but the methods apply more generally.

(2)  The assumption $n<m+\frac{n-1}{2}$ is unnecessarily strong; extension down to the center of symmetry of the functional equation should be possible by the methods of \cite{H08}.

(3) It is likely that this result can be improved to allow for the action of $Gal(\Qbar/F')$ on the ratio of the two sides, for some subfield $F' \subset F^{gal}$; one would like to replace $F^{gal}$ by $\QQ$.

\end{rmks}

\subsection{Results on Rankin-Selberg $L$-functions}
Let $\Pi$ be a cohomological regular polarized cuspidal representation of $GL_{n}(\AF)$. We assume that $\Pi^{\vee}$ descends to $\{\pi\}_V$, a packet of representations of the rational similitude group associated to $V$, which contains a holomorphic discrete series representation and we denote it $\pi$. The automorphic period $Q_{V}(\pi)$ can be defined as before.
\begin{df}\label{autoperiods}
Let $I=(r_{\sigma})_{\sigma\in\Sigma}$ be an element in $\{0,1,\cdots,n\}^{\Sigma}$. We define the \textbf{automorphic arithmetic period} $P^{(I)}(\Pi):=(2\pi)^{-2a_{0}}Q_{V}(\pi)$.

If $F^+ = \QQ$ then we write $P^{(a)}(\Pi)$ for $P^{(I)}(\Pi)$ when $I$ is the singleton $a$.
\end{df}

We observe that \begin{equation}\nonumber
L\left(m,\pi\otimes\psi, St\right) = L\left(m,\Pi^{\vee}\otimes \widetilde{\psi}\right) = L\left(m,\Pi^{c}\otimes \widetilde{\psi}\right) = L\left(m,\Pi\otimes \widetilde{\psi^{c}}\right).
\end{equation}

If we have moreover that $\Pi^{c}\cong \Pi$ then $\pi$ is polarized and we may apply Guerberoff's result to the left hand side. In this case, we have
\begin{thm}
Let $\Pi$ be as the beginning of this subsection. If $\Pi$ moreover satisfies $\Pi^{c}\cong \Pi$, then for $ m\in \Z+\frac{n-1}{2}$ satisfying (\ref{critical condition gue}), we have:
\begin{eqnarray}
& L^{S}\left(m,\Pi\otimes \widetilde{\psi^{c}}\right)\sim_{E(\Pi)\otimes E(\psi);F^{gal}}\nonumber&\\
&(2\pi i)^{d(F^{+})mn/2}P^{(I)}(\Pi)\prod\limits_{\sigma\in\Sigma}Q^{(0)}(M(\widetilde{\psi^{c}}),\sigma)^{n-I_{\sigma}}Q^{(1)}(M(\widetilde{\psi^{c}}),\sigma)^{I_{\sigma}}&
\end{eqnarray}
where $I_{\sigma}=r_{\sigma}$ in (\ref{critical condition gue}).
\end{thm}

We remark that $E(\pi) \supset E(\Pi)$ but the inclusion is in general strict. 
But since all terms on both sides other than $P^{(I)}(\Pi)$ depend only on the embedding of $E(\Pi)\otimes E(\psi)$ into $\C$, we know that the arithmetic automorphic period $P^{(I)}(\Pi)$ also depends only on the embedding of $E(\Pi)\otimes E(\psi)$. Moreover, we may replace $E(\pi)\otimes E(\psi)$ by $E(\Pi)\otimes E(\psi)$ in the relation.

If we arrange $(a_{i,\sigma})_{1\leq i\leq n}$ in decreasing order for any fixed $\sigma\in\Sigma$, then the infinity type of $\Pi$ at $\sigma\in\Sigma$ is $z^{-a_{i,\sigma}-\frac{n+1}{2}+i}\overline{z}^{a_{i,\sigma}+\frac{n+1}{2}-i}$, and the infinity type of $\widetilde{\psi^{c}}$ is $z^{m_{\sigma}-m_{\overline{\sigma}}}\overline{z}^{- m_{\sigma}+m_{\overline{\sigma}}}$. It is easy to see that equation (\ref{critical condition gue}) implies that $r_{\sigma}=I_{\sigma}(\Pi, \widetilde{\psi^{c}})$ defined in Proposition \ref{n1motive}. 

The methods of the above theorem should work in more general cases. As mentioned above, we expect that the field $F^{gal}$ can be replaced by $\Q$.  We assume it as a hypothesis here. L. Guerberoff and the second author plan to include a proof in a future paper.

\begin{hyp}\label{hypCrelle}
The above theorem is true for general $\Pi$ as in the beginning of this subsection and all critical $m>0$ with $I_{\sigma}=I_{\sigma}(\Pi, \widetilde{\psi^{c}})$ and the relation up to $E(\Pi)\otimes E(\psi)$.
\end{hyp}

This hypothesis is the Conjecture $5.1.1$ of \cite{linthesis}. Under this hypothesis, one can prove the following factorization result as in Theorem $7.6.1$ of \textit{loc.cit.}

\begin{thm}\label{factorizationtheorem}
If $\Pi$ is sufficiently regular then there exist some complex numbers $P^{(r)}(\Pi,\sigma)$, $0\leq r\leq n$, well-defined up to multiplication by elements in $(E(\Pi)\sigma(F))^{\times}$, such that the following two conditions are satisfied:
\begin{enumerate}
\item $P^{(I)}(\Pi) \sim_{E(\Pi); F^{gal}} \prod\limits_{\sigma\in\Sigma}P^{(I(\sigma))}(\Pi,\sigma)$ for all $I=(I(\sigma))_{\sigma\in\Sigma}\in \{0,1,\cdots,n\}^{\Sigma}$,
\item  and $P^{(0)}(\Pi,\sigma)\sim_{E(\Pi); F^{gal}} p(\widecheck{\xi_{\Pi}},\overline{\sigma})$
\end{enumerate}
where $\xi_{\Pi}$ is the central character of $\Pi$.

Moreover, we know $P^{(n)}(\Pi,\sigma)\sim_{E(\Pi); F^{gal}} p(\widecheck{\xi_{\Pi}},\sigma)$ or equivalently $P^{(0)}(\Pi,\sigma)\times P^{(n)}(\Pi,\sigma)\sim_{E(\Pi); F^{gal}} 1$.
\end{thm}

The following two theorems are known if $F$ is a quadratic imaginary field. The generalization to CM fields is immediate once we know Hypothesis \ref{hypCrelle} and the above factorization theorem. We refer to section $9.5$ of \cite{linthesis} for a discussion of the generalization. We recall that the proof of Theorem \ref{factorizationtheorem} is also based on Hypothesis \ref{hypCrelle}. 

 \begin{thm}\label{RSGood}  
Let $\Pi$ and $\Pi'$ be cuspidal cohomological automorphic representation of $GL_{n}(\AF)$ and $GL_{n'}(\AF)$ respectively which satisfies the descending condition as in the beginning of Section $4.3$. We assume that both $\Pi$ and $\Pi'$ are sufficiently regular and that $(\Pi,\Pi')$ is in good position. If $n \nequiv n' \mod 2$, then for positive $m\in \Z+\frac{n+n'}{2}$ which is critical for $\Pi\times \Pi'$, we have
\begin{eqnarray}
&L(m,\Pi\times \Pi') \sim_{E(\Pi)\otimes E(\Pi'); F^{gal}}&\\
 &(2\pi i)^{mnn'd(F^{+})}\prod\limits_{\sigma\in\Sigma} [\prod\limits_{j=0}^{n}(P^{(j)}(\Pi,\sigma))^{sp(j,\Pi;\Pi',\sigma)}\prod\limits_{k=0}^{n'}(P^{(k)}(\Pi',\sigma))^{sp(k,\Pi';\Pi,\sigma)}].&\nonumber
\end{eqnarray}

\end{thm}

For more general configurations, we can show the following result.

\begin{thm}\label{RSgeneral}  
Theorem \ref{RSGood} is still true without the good position condition for $m=1$ when $n\equiv n' \mod 2$.
\end{thm}

\begin{rem}
We have similar results for general relative parity of n and n'. We refer to Theorem $10.8.1$ (resp. Theorem $11.4.1$) of \cite{linthesis} for the precise statement of Theorem \ref{RSGood} (resp. Theorem \ref{RSgeneral}).
\end{rem}

\subsection{Geometric meaning of local periods}

Let $r$ be an integer between $0$ and $n$. Let $I\in \{0,1,\cdots,n\}^{\Sigma}$ such that $I(\sigma)=r$ and $I(\sigma')=n$ for all $\sigma'\neq \sigma$. By Theorem \ref{factorizationtheorem} we have 
\begin{equation}\label{decomposition of automorphic period}
P^{(I)}(\Pi) \sim_{E(\Pi); F^{gal}} P^{(r)}(\Pi,\sigma)\prod\limits_{\sigma'\neq \sigma} P^{(n)}(\Pi,\sigma') \sim_{E(\Pi); F^{gal}} P^{(r)}(\Pi,\sigma)\prod\limits_{\sigma'\neq \sigma} p(\widecheck{\xi_{\Pi}},\sigma').
\end{equation}

We recall that $P^{(I)}(\Pi)$ relates to the representation of the similitude unitary group with base change $\Pi^{\vee}\otimes \xi$ where $\xi$ is a Hecke character over $F$ such that $\xi_{\Pi}=\cfrac{\xi}{\xi^{c}}$. Therefore, $P^{(I)}(\Pi)$ should be equivalent to the inner product of a rational class in the bottom stage of $\Lambda^{n-r}M_{\sigma}(\Pi^{\vee}) \otimes M(\xi)$. More precisely, we should have
\begin{eqnarray}
&&P^{(I)}(\Pi) \nonumber\\
&\sim_{E(\Pi)\otimes E(\xi); F^{gal}} &Q_{1}(M(\Pi^{\vee}),\sigma)Q_{2}(M(\Pi^{\vee}),\sigma)\cdots Q_{n-r}(M(\Pi^{\vee}),\sigma) Q_{1}(M(\xi))\nonumber
\\
&\sim_{E(\Pi)\otimes E(\xi); F^{gal}} &Q_{1}(M(\Pi^{\vee}),\sigma)Q_{2}(M(\Pi^{\vee}),\sigma)\cdots Q_{n-r}(M(\Pi^{\vee}),\sigma) Q_{1}(M(\xi),\sigma) \prod\limits_{\sigma'\neq \sigma} Q_{1}(M(\xi),\sigma')\nonumber
\end{eqnarray}

By Lemma \ref{comparison n=1}, we know 
\begin{equation}\label{a calculation for sigma'}
Q_{1}(M(\xi),\sigma') \sim _{E(\xi)} p( \cfrac{\widecheck{\xi}}{\widecheck{\xi^{c}}},\sigma') \sim _{E(\xi)}  p(\widecheck{\xi_{\Pi}},\sigma'). 
\end{equation}

We compare this with equation (\ref{decomposition of automorphic period}) and then deduce that the Tate conjecture would imply (as in the section $4.5$ of \cite{GH}):
\begin{eqnarray}\label{main comparison 1}
&&P^{(r)}(\Pi,\sigma) \\\nonumber
&\sim_{E(\Pi); F^{gal}}& Q_{1}(M(\Pi^{\vee}),\sigma)Q_{2}(M(\Pi^{\vee}),\sigma)\cdots Q_{n-r}(M(\Pi^{\vee}),\sigma) Q_{1}(M(\xi),\sigma)\\
\nonumber
&\sim_{E(\Pi); F^{gal}}& Q_{1}(M(\Pi^{c}),\sigma)Q_{2}(M(\Pi^{c}),\sigma)\cdots Q_{n-r}(M(\Pi^{c}),\sigma) Q_{1}(M(\xi),\sigma)
\end{eqnarray}

We repeat the calculation in (\ref{a calculation for sigma'}) for $\sigma$ and get $Q_{1}(M(\xi),\sigma) \sim_{E(\xi)} p(\widecheck{\xi_{\Pi}},\sigma) \sim_{E(\xi)}  \delta(M(\xi_{\Pi}^{c}),\sigma)$ by Lemma \ref{comparison n=1}.

The fact that $det(M(\Pi^{c})) \cong det(M(\xi_{\Pi^{c}}))(\frac{n(1-n)}{2})$ implies $\delta(M(\Pi),\sigma) \sim_{E(\Pi); F^{gal}} \delta(M(\xi_{\Pi}^{c}),\sigma) (2\pi i)^{n(1-n)/2}$ and hence $Q_{1}(M(\xi),\sigma)\sim_{E(\Pi); F^{gal}} \delta(M(\Pi),\sigma)(2\pi i)^{n(n-1)/2}$.

This implies that $P^{(r)}(\Pi,\sigma)\sim_{E(\Pi); F^{gal}} Q^{(n-r)}(M(\Pi^{c}),\sigma)$ where the right hand side is equivalent to $Q^{(r)}(M(\Pi),\sigma)$ by Lemma \ref{conjugacy on motivic period} as expected.

\section{Integral representations}

Let $G = GL(n)$, $G' = GL(n-1)$ over the number field $F$, with notation as in the introduction.  As above, $F$ is assumed to be a CM field.
Let $\Pi \times \Pi'$ be an automorphic representation of $G \times G'$ over $F$, with $\Pi$ cuspidal, and let
$$\iota:  G' \hookrightarrow G,  ~\iota(g') = diag(g',1).$$
Write $\Pi = \Pi_{\infty}\otimes \Pi_f$, and likewise for $\Pi'$, the factorization of the adelic representations into their archimedean and non-archimedean components.
 
 Let $\phi \in \Pi$, $\phi' \in \Pi'$.
We normalize the Jacquet-Piatetski-Shapiro-Shalika zeta integral for $\Pi \times \Pi'$ over $F$: 
$$Z(s,\phi,\phi') = \int_{G'(F)\backslash G'(\ad)} \phi(\iota(g'))\phi'(g') ||det(g')||^{s - \frac{2n-3}{2}} dg'.$$
With this normalization, the central point of the functional equation is at $s = \frac{2n-3}{2}$.

The integral admits an Euler product factorization, in the usual way.  Suppose $\phi = \otimes_v \phi_v$ and $\phi' = \otimes_v \phi'_v$ with respect to given factorizations
$\Pi \isoarrow \otimes'_v \Pi_v$;  $\Pi' \isoarrow \otimes'_v \Pi'_v$.  Then
\begin{equation}\label{euler} Z(s,\phi,\phi') = \prod_v Z_v(s,\phi_v,\phi'_v). \end{equation}

When $\Pi$ and $\Pi'$ are cohomological representations, and $\phi$, $\phi'$ cohomological vectors,
$Z(s,\phi,\phi')$ can be interpreted as a cup product on the locally symmetric space for $G \times G'$.   Suppose 

\begin{itemize}
\item [(a)] $\Pi$ is cuspidal;
\item  [(b)]  $\Pi'_{\infty}$ is tempered (up to a twist by a power of the determinant) and generic;
\item  [(c)]  $H^i(\mathfrak{gl}(n),U(n); \Pi_{\infty}\otimes W) \neq 0$, $H^{i'}(\mathfrak{gl}(n-1),U(n-1); \Pi'_{\infty}\otimes W')\neq 0$ for some irreducible finite dimensional representations $W$ and $W'$ of $R_{F/\QQ}G$ and  $R_{F/\QQ}G'$, respectively and for some $i$ and $i'\in\Z$.
\end{itemize}

There is an interval $I(n,F) = [b_n(F),t_n(F)] \subset \ZZ$ such that 
$$H^i(\mathfrak{gl}(n),U(n); \Pi_{\infty}\otimes W) \neq 0 \text{ if and only if } i \in I(n,F).$$ 
If $F$ is a CM field of degree $2d$ over $\QQ$, then $b_n(F) = \frac{n(n-1)}{2}\cdot d$ and  $\dim H^{b_n} = 1$.   By hypothesis (b) above, the same holds for $G'$ and $W'$.   

Note that the intervals $I(n,F)$ and $I(n-1,F)$ do not depend on the coefficients $W, W'$.  We let $E(W)$ and $E(W')$ denote the fields of definition of $W$ and $W'$ as representations of $R_{F/\QQ}G$ and  $R_{F/\QQ}G'$, respectively; these are always finite extensions of $\QQ$.  Let $E(W,W') = E(W)\cdot E(W')$.

In general, $b_n(F) + b_{n-1}(F) = \dim (R_{F/\QQ}G')(\RR)/U(n-1)^d = \dim {}_K \tilde{S}_{n-1}$, for any level subgroup $K \subset G'(\ad^f)$, where
$${}_K \tilde{S}_{n-1} = G'(F)\backslash R_{F/\QQ}G'(\ad)/U(n-1)^d\times K, ~  K \subset G'(\ad^f).$$
Note that this is {\it not} the locally symmetric space attached to $R_{F/\QQ}G'$, because we have not taken the quotient by the center of $R_{F/\QQ}G'(\RR)$.
Nevertheless, the expression $\phi(\iota(g'))\phi'(g')$ can be identified with a top degree differential, that we may denote $\omega_{\phi} \cup \omega_{\phi'}$ on ${}_K \tilde{S}_{n-1}$
for appropriate choices of $\phi, \phi'$.  The JPSS integral is obtained by taking the image of this differential in compactly-supported cohomology in the top degree and pairing with the Borel-Moore homology class defined by ${}_K \tilde{S}_{n-1}$ itself.  The pairing is well-defined over the field of definition $E(W,W')$.    The action of Hecke operators preserve the $E(W)$ and $E(W')$-rational structures  of the cohomology, and by multiplicity one for $GL(n)$, the subspaces of cohomology defined by $\Pi_f$ and $\Pi'_f$ are defined over finite extensions $E(\Pi) \supset E(W)$, $E(\Pi') \supset E(W')$.  

In particular, suppose $\phi$ and $\phi'$ are chosen so that $\omega_{\phi}$ and $\omega_{\phi'}$ define $E(\Pi)$ and $E(\Pi')$-rational cohomology classes, respectively.  Suppose moreover that there is a non-zero $R_{F/\QQ}G'$-equivariant homomorphism
\begin{equation}\label{contraction}  \xi:  W \otimes W' \ra Triv, \end{equation}
where $Triv$ is the trivial one-dimensional representation of $R_{F/\QQ}G'$; we may assume that $\xi$ is rational over $E(W,W')$.   The cup product
$\omega_{\phi} \cup \omega_{\phi'}$ naturally belongs to 
$$H^{b_n(F) + b_{n-1}(F)}({}_K \tilde{S}_{n-1},\iota^*(\tilde{W})\otimes \tilde{W'})$$
where $\tilde{W}$ and $\tilde{W'}$ are the local systems on the locally symmetric space attached to $G$ and on ${}_K \tilde{S}_{n-1}$, respectively, attached to the representations $W$ and $W'$, and $\iota^*(\tilde{W})$ is the pullback of $\tilde{W}$ to ${}_K \tilde{S}_{n-1}$ by the map defined by $\iota:  G' \hookrightarrow G$.  Applying $\xi$, we find that $\xi(\omega_{\phi} \cup \omega_{\phi'}) \in H^{\dim {}_K \tilde{S}_{n-1}}({}_K \tilde{S}_{n-1},Triv)$, in the obvious notation.
Then the integral at the central point $s = \frac{2n-3}{2}$ is exactly the pairing of the cup product $\omega_{\phi} \cup \omega_{\phi'}$ with the top class in Borel-Moore homology of ${}_K \tilde{S}_{n-1}$, and thus is an element of $E(\Pi,\Pi') = E(\Pi)\cdot E(\Pi')$.   

More generally, let $W(m) = W\otimes \det^{-m}$.  If $m$ is an integer such that there are equivariant contractions
\begin{equation}\label{contractionm}  \xi(m):  W(m) \otimes W' \ra Triv,   \xi^{\vee}(m):  W^{\vee}(m) \otimes W^{\prime,\vee}   \ra Triv, \end{equation}
where the superscript $^{\vee}$ denotes contragredient, then the JPSS integral at $s$ can also be interpreted as a cohomological cup product, and thus is again in $E(\Pi,\Pi')$.

Condition \eqref{contraction} corresponds to the Good Position Hypothesis.  When $F$ is imaginary quadratic, the set of $m$ satisfying the property \eqref{contractionm} is identified in Lemma 3.5 of \cite{GH} with the set of critical points of the $L$-function $L(s,\Pi \otimes \Pi')$ to the right of the center of symmetry of the functional equation (including the center if it is an integer).   The same calculation holds for general CM fields.

The rationality property of the cup product does not respect the Euler product factorization.  The right-hand side of \eqref{euler} is a product of local integrals defined in terms of the Whittaker models of the local components $\Pi_v$, $\Pi'_v$.  In order to compare the local integrals -- especially the local factors at unramified places --   to standard Euler factors, that depend only on the local representations, one needs to introduce {\it Whittaker periods} 
$$p(\Pi) \in (E(\Pi)\otimes \CC)^{\times}, ~~ p(\Pi')  \in (E(\Pi')\otimes \CC)^{\times}$$ 
that measure the difference between the rational structure defined by cohomology and that defined by standard Whittaker models.  There is a single natural definition of the latter at non-archimedean places; at archimedean places one makes an arbitrary choice.  In the end, one obtains a formula of the form
\begin{equation}\label{comp1} Z_S(m,\phi,\phi')L(m,\Pi,\Pi') = Z(m,\phi,\phi') \in p(\Pi)p(\Pi')p(m,\Pi_{\infty},\Pi'_{\infty})\cdot E(\Pi,\Pi') \end{equation}
when $m$ satisfies \eqref{contractionm}. 
Manipulating this expression, one obtains a preliminary version of Theorem \ref{RSoverQ}:
\begin{equation}\label{comp2} L(m,\Pi,\Pi') \sim_{E(\Pi)E(\Pi')}  p(\Pi)p(\Pi')Z(m,\Pi_{\infty},\Pi'_{\infty}). \end{equation}
(The last two formulas correspond roughly to Theorem 3.9 of \cite{GH}.)

The next step is to express the Whittaker periods in terms of the automorphic motivic periods $P^{(I)}(\Pi)$.  This is done in \cite{GH} and \cite{linthesis} in several steps, based on choosing $\Pi'$ in appropriate spaces of Eisenstein cohomology.  Fix a cuspidal cohomological $\Pi$ and suppose for simplicity that $n$ is even.  One can find a cohomological $\Pi'$ consisting of Eisenstein classes attached to a Hecke character $\chi = (\chi_1,\dots,\chi_{n-1})$ of the Levi subgroup $GL(1)^{n-1}$ of a Borel subgroup of $G'$.  Then we have
\begin{equation}\label{piprime} L(s,\Pi,\Pi') = \prod_{i = 1}^{n-1}  L(s,\Pi,\chi_i).\end{equation}
If the $\chi_i$ are chosen appropriately, $\Pi'$ is cohomological with respect to a $W'$ that satisfies \eqref{contraction}; thus $L(m,\Pi,\Pi')$ for critical $m$ can be related to $p(\Pi)\cdot p(\Pi')$.  On the other hand, Theorem \ref{guer} expresses the critical values of the right-hand side of \eqref{piprime} in terms of the motivic automorphic periods $Q_{V_i}(\Pi)$ of $\Pi$, for $V_i$ of varying signature, and of the Hecke characters $\chi_i$.  Finally, Shahidi's formulas for the Whittaker coefficients of Eisenstein series express $p(\Pi')$ in terms of the same automorphic periods of the $\chi_i$.   In the end, one finds that
\begin{equation}\label{whittakerperiods}  p(\Pi) \sim_{E(\Pi); F^{gal}} (*) \prod_{i = 1}^{n-1}Q_{V_i}(\Pi) \end{equation}
where $(*)$ is an elementary factor.  

Now assume $\Pi'$ is cuspidal.  Combining \eqref{whittakerperiods} with \eqref{comp2}, one obtains an expression for the critical values under the Good Position Hypothesis in terms of motivic automorphic periods of the form $Q_{V_i}(\Pi)$  and $Q_{V'_j}(\Pi')$.  Using period relations for automorphic induction, Lin was able in \cite{LinCR} to replace the unspecified elementary and archimedean factors by explicit powers of $2 \pi i$.  This proves the automorphic version \ref{RSoverQ} of Conjecture \ref{deligneMMprime} in the Good Position situation, and its extension \ref{RSGood} to general CM fields. 

The proof of the automorphic theorem \ref{RSgeneral} in \cite{linthesis} is based on the same principle, except this time the interpretation involves a comparison of an unknown factor in Shahidi's formula for the Whittaker coefficient with a motivic expression.

Finally, the proof in \cite{linthesis} of Theorem \ref{factorizationtheorem} is based on comparing several expressions of the form \eqref{comp1} for critical values, that are obtained by interpreting Rankin-Selberg integrals as cohomological cup products.  The regularity hypothesis in Theorem \ref{factorizationtheorem} is required to guarantee that none of the $L$-functions used in the comparison vanishes at the relevant points.

\appendix

\begin{section}{Tensor product and Compositions}
Throughout the text, fix $E$ a number field. 
Let $F$ be a field containing $\Q$. In the applications, $F$ will be either a number field or the complex field $\C$.

We denote by $\Sigma_{E}$ (resp. $\Sigma_{F}$) the set of embeddings of $E$ (resp. $F$) in $\C$.

For $\sigma\in\Sigma_{F}$, we write $\overline{\sigma}$ for the complex conjugation of $\sigma$. 

Tensor products without subscript are by default over $\Q$. We also write $\otimes_{\sigma}$ for $\otimes_{F,\sigma}$ etc. For example, if $V$ is an $F$-vector space, we write $V\otimes_{\sigma}\C $ for $V\otimes_{F,\sigma}\C $.

\begin{df}
Let $(L,\iota,\nu)$ be a triple where $L$ is a field, $\iota$ is an embedding of $E$ in $L$ and $\nu$ is an embedding of $F$ in $L$. We say this triple is a \textbf{compositum} of $E$ and $F$ if $\iota(E)$ and $\nu(F)$ generate $L$.

Two compositums $(L,\iota,\nu)$ and  $(L',\iota',\nu')$ are called isomorphic if there exists an isomorphism of fields $L\cong L'$ which commutes with the embeddings.
\end{df}

\begin{prop}\label{decomposition algebra}
\begin{enumerate}
\item The $\Q$-algebra $E\otimes F$ decomposes uniquely as a direct sum of fields:
\begin{equation}
E\otimes F =\bigoplus\limits_{\alpha \in \mathcal{A}} L_{\alpha}
\end{equation}
\item For each $\alpha\in \mathcal{A}$, let $\iota_{\alpha}$ (resp. $\nu_{\alpha}$) be the composition of the canonical map $E\rightarrow E\otimes F$ (resp. $F\rightarrow E\otimes F$) and the projection of $E\otimes F$ to $L_{\alpha}$. The triple $(L_{\alpha},\iota_{\alpha},\mu_{\alpha})$ is a compositum of $E$ and $F$.
\item If $\alpha$ and $\alpha'$ are two different elements in $\mathcal{A}$, then the two compositums $(L_{\alpha},\iota_{\alpha},\mu_{\alpha})$ and $(L_{\alpha'},\iota_{\alpha'},\mu_{\alpha'})$ are not isomorphic.
\item Any compositum of $E$ and $F$ can be obtained in this way up to isomorphisms.
\end{enumerate}
\end{prop}
\begin{dem}
\begin{enumerate}
\item The uniqueness is clear. Let us prove the existence.

Write $E\cong \Q[X]/(f)$ where $f$ is an irreducible polynomial in $\Q[X]$. Then the $\Q$-algebra $E\otimes F \cong F[X]/(f)$. 

We decompose $f=\prod\limits_{\alpha \in \mathcal{A}}f_{\alpha}$ in the ring $F[X]$. Since $f$ is separable we know that the polynomials $f_{\alpha}$, $\alpha\in\mathcal{A}$, are different.

Therefore, $F[X]/(f)\cong \bigoplus\limits_{\alpha\in \mathcal{A}} L_{\alpha}$ where $L_{\alpha}=F[X]/(f_{\alpha})$ is a field as predicted.

\item The map $\iota_{\alpha}$ is a map from $\Q[X]/(f)$ to $F[X]/(f_{\alpha})$ sending $X$ to $X$. The map $\nu_{\alpha}:F\rightarrow F[X]/(f_{\alpha})$ is induced by the natural embedding $F\rightarrow F[X]$. It is easy to see that the image of $\iota_{\alpha}$ and $\nu_{\alpha}$ generate $F[X]/(f_{\alpha})=L_{\alpha}$.
\item This is due to the fact that the polynomials $f_{\alpha}$, $\alpha\in\mathcal{A}$, are different.
\item Let $(L,\iota,\nu)$ be a compositum of $E$ and $F$. Then the map $\iota\otimes\nu:E\otimes F\rightarrow L$ induces a surjective ring homomorphism from $\bigoplus\limits_{\alpha \in \mathcal{A}} L_{\alpha}$ to $L$. Since $L$ is a field, this ring homomorphism must factor through one of the $L_{\alpha}$ and then (4) follows.

\end{enumerate}

\end{dem}

\begin{rem}
 It is easy to see that there is a bijection between $\Sigma_{E}\times \Sigma_{F}$ and $\bigsqcup\limits_{\alpha\in\mathcal{A}} \Sigma_{L_{\alpha}}$. More precisely, let $\tau\in\Sigma_{E}$ and $\sigma\in\Sigma_{F}$, then $\tau(E)$ and $\sigma(F)$ generate a number field $L$. The triple $(L,\tau,\sigma)$ is a compositum of $E$ and $F$ and hence is isomorphic to $(L_{\alpha},\iota_{\alpha},\nu_{\alpha})$ for a unique $\alpha\in\mathcal{A}$. We define $\alpha(\tau,\sigma)=\alpha$.
\end{rem}

\begin{rem}\label{index is Galois invariant}
For any $g\in Gal(\bar{\Q}/\Q)$, we have $\alpha(g\tau,g\sigma)=\alpha(\tau,\sigma)$. In fact, let $L$ be the field generated by $\tau(E)$ and $\sigma(F)$. Then the compositum $(L,\tau,\sigma)$ is isomorphic to $(gL,g\tau,g\sigma)$. By the third point of Proposition \ref{decomposition algebra}, we know $\alpha(g\tau,g\sigma)=\alpha(\tau,\sigma)$.

On the other hand, it is easy to see that for any $\sigma,\sigma'\in \Sigma_{F}$ and $\tau,\tau'\in\Sigma_{E}$, $\alpha(\tau,\sigma)=\alpha(\tau',\sigma')$ implies that there exists $g\in Gal(\bar{\Q}/\Q)$ such that $(g\tau,g\sigma)=(\tau',\sigma')$. 

In other words, the isomorphism classes of compositums are in bijection with the $Gal(\bar{\Q}/\Q)$-orbits of $\Sigma_{E}\times \Sigma_{F}$.

\end{rem}

It is easy to prove the following lemma:
\begin{lem}\label{composition lemma}
\begin{enumerate}
Let $\tau\in\Sigma_{E}$ and $\sigma\in\Sigma_{F}$.
\item The $\C$-vector space $L_{\alpha}\otimes_{\tau\otimes \sigma} \C\neq 0$ if and only if $\alpha=\alpha(\tau,\sigma)$. 
\item The following equation holds:
\begin{equation}
L_{\alpha(\tau,\sigma)}\otimes_{\tau\otimes \sigma} \C=(E\otimes F)\otimes_{\tau\otimes \sigma} \C
\end{equation}
\end{enumerate}
\end{lem}

\end{section}

\begin{section}{Decomposition of $E\otimes F$-modules}
The decomposition of $E\otimes F$ as $\Q$-algebra in Proposition \ref{decomposition algebra} gives a decomposition on $E\otimes F$-modules as follows:
\begin{prop}
Let $M$ be an $E\otimes F$-module. It decomposes as a direct sum
\begin{equation} \label{decomposition via alpha}
M=\bigoplus\limits_{\alpha \in \mathcal{A}} M(\alpha)
\end{equation}
where $M(\alpha)$ is an $L_{\alpha}$-vector space and the action of $E\otimes F$ on $M(\alpha)$ factors through the action of $L_{\alpha}$.
\end{prop}

Lemma \ref{composition lemma} then implies that:
\begin{prop} \label{prop for decomposition}
Let $M$ be a finitely generated $E\otimes F$-module, $\sigma\in\Sigma_{F}$ and $\tau\in \Sigma_{E}$.
\begin{enumerate}
\item If $\alpha\neq \alpha(\tau,\sigma)$ is an element in $\mathcal{A}$ then $M(\alpha)\otimes_{\tau\otimes \sigma} \C = 0$. 
\item If $\alpha= \alpha(\tau,\sigma)$ then $M\otimes_{\tau\otimes \sigma} \C = M(\alpha)\otimes_{\tau\otimes \sigma} \C $. 
\item The dimension of $M\otimes _{\tau\otimes \sigma} \C$ over $\C$ is equal to $dim_{L_{\alpha(\tau,\sigma)}}M(\alpha(\tau,\sigma))$. In particular, it only depends on $\alpha(\tau,\sigma)$.
\item The module $M$ is free if and only if $dim_{\C} M\otimes _{\tau\otimes \sigma} \C$ is the same for all $(\tau,\sigma)$.
\end{enumerate}
\end{prop}

This lemma can be easily deduced from the fact that $M(\alpha)\otimes_{\tau\otimes \sigma} \C=M(\alpha)\otimes_{F_{\alpha}} F_{\alpha}\otimes_{\tau\otimes \sigma} \C$  for any $\alpha\in \mathcal{A}$.

\begin{lem}\label{basislemma}
Let $M$ be a free $E\otimes F$-module of rank $n$. A family of $n$ element in $M$ forms an $E\otimes F$ basis is equivalent to that it is linearly independent over $E\otimes F$, and is also equivalent to that it generates $M$ over $E\otimes F$.
\end{lem}
This can be deduced from similar results for vector spaces over the fields $L_{\alpha}$.

\bigskip

The previous results, when applied to the case $F=\C$, give a decomposition for any $E\otimes \C$-module. More precisely, we identify the $\Q$-algebra $E\otimes \C$ with $\C^{\Sigma_{E}}$ by sending $e\otimes z$ to $(\tau(e)z)_{\tau\in\Sigma_{E}}$ for all $\tau \in E$ and $z\in \C$. We then have:

\begin{prop}
Let $V$ be a finitely generated $E\otimes \C$-module. The action of $E$ gives a decomposition of $V$ as direct sum of sub-$E\otimes \C$-modules:
\begin{equation}
V=\bigoplus\limits_{\tau:E\hookrightarrow \C}V(\tau)
\end{equation}
where the action of E on $V(\tau)$ is given by scalar multiplication via $\tau$.
\end{prop}
For example, if $M$ is an $E$ vector space, then the $E\otimes \C$-module $M \otimes \C =\bigoplus\limits_{\tau:E\hookrightarrow \C} M\otimes _{\tau} \C$.

\begin{prop}\label{basis component}
\begin{enumerate}
\item A finitely generated $E\otimes \C$-module $V$ is free if and only if $dim_{\C}V(\tau)$ is the same for all $\tau$. 

\item Let $V$ be a free finitely generated $E\otimes \C$-module of rank $d$. By $(1)$, each $V(\tau)$ has dimension $d$. For each $\tau\in \Sigma_{E}$, let $\{w_{1}(\tau),w_{2}(\tau),\cdots, w_{d}(\tau)\}$ be a $\C$-basis of $V(\tau)$. For each $1\leq i\leq d$, we put $w_{i}=\sum\limits_{\tau\in \Sigma_{E}}w_{i}(\tau)\in V$. Then the family $\{w_{1},w_{2},\cdots,w_{d}\}$ forms an $E\otimes \C$-basis of $V$.
\item Let $V_{1}$ and $V_{2}$ be two free finitely generated $E\otimes \C$-modules. Let $f:V_{1} \rightarrow V_{2}$ be an isomorphism of $E\otimes \C$-modules. For each $\tau$, it induces an isomorphism $f(\tau):V_{1}(\tau) \rightarrow V_{2}(\tau)$.

 If for each $\tau$, we have fixed $\C$-bases of $W_{1}(\tau)$ and $W_{2}(\tau)$, we can construct $E\otimes \C$-bases for $W_{1}$ and $W_{2}$ as in $(2)$. With respect to these fixed bases, we have: $det(f)=(det(f(\tau)))_{\tau\in \Sigma_{E}}\in E\otimes \C$.

\end{enumerate}
\end{prop}

\end{section}

\end{document}